\documentclass[12pt,oneside,a4paper,reqno]{amsart}

\usepackage{amssymb}
\usepackage{amsmath}
\usepackage{amsthm}
\usepackage{amscd}
\usepackage[all]{xy}
\usepackage{longtable}
\usepackage{mathrsfs}

\usepackage{bm}

\usepackage{comment}

\usepackage{graphicx,xcolor}
\usepackage{pict2e}

\usepackage{tikz}
\usetikzlibrary{cd}

\usepackage{todonotes}

\usepackage{hyperref}
\hypersetup{
    colorlinks=true,
    linkcolor=red,
    citecolor=blue
}

\usepackage{tikz}
\usetikzlibrary{cd}

\usepackage{geometry}
\theoremstyle{plain}
\newtheorem{theorem}{Theorem}[section]
\newtheorem{lemma}[theorem]{Lemma}
\newtheorem{corollary}[theorem]{Corollary}
\newtheorem{proposition}[theorem]{Proposition}
\newtheorem{remark}[theorem]{Remark}

\newtheorem*{mainthm}{Main Theorem}

\makeatletter
\newcommand{\mainthmlabel}[1]{%
  \phantomsection
  \def\@currentlabel{Main Theorem}%
  \label{#1}%
}
\makeatother

\theoremstyle{definition}
\newtheorem{definition}[theorem]{Definition}

\newtheorem{remark-theorem}[theorem]{Remark-Theorem}

\newtheorem{notation}[theorem]{Notation}

\newcounter{setting}

\newcommand{\kah}{K\"{a}hler }
\newcommand{\idd}{i\partial\overline{\partial}}
\newcommand{\dbar}{\overline{\partial}}
\newcommand{\delh}{\partial_h}
\newcommand{\Ker}{\rom{Ker}\,}
\newcommand{\Dom}{\rom{Dom}\,}
\newcommand{\Imege}{\rom{Im}\,}

\newcommand{\cal}[1]{\mathcal{#1}}
\newcommand{\bb}[1]{\mathbb{#1}}
\newcommand{\scr}[1]{\mathscr{#1}}
\newcommand{\rom}[1]{\mathrm{#1}}

\newcommand{\tx}{\widetilde{X}}

\newcommand{\tl}[1]{\widetilde{#1}}
\newcommand{\wh}[1]{\widehat{#1}}

\newcommand{\calH}{\cal{H}^{n+m,q}(X,L,h,\Phi)}
\newcommand{\lih}{L\otimes\scr{I}(h)}
\newcommand{\liho}{L_0\otimes\scr{I}(h_0)}
\newcommand{\hdi}{R^qf_*(K_{X/Y}\otimes\lih)}

\newcommand{\calHo}{\cal{H}^{n,q}(X_0,L_0,h_0)}
\newcommand{\Lxo}[1]{L^2_{n,{#1}}(X_0,L_0,h_0)}

\newcommand{\hajo}{\widehat{\alpha}^j_0}
\newcommand{\hajoe}{\widehat{\alpha}_{0,\varepsilon}^j}
\newcommand{\bjo}{\beta^j_0}
\newcommand{\cjo}{\gamma^j_0}
\newcommand{\ujo}{u^j_0}

\newcommand{\Tqe}{\dbar_{q,\varepsilon}}
\newcommand{\cel}{\gamma_{0,\varepsilon,\lambda}^j}
\newcommand{\wel}{w_{0,\varepsilon,\lambda}^j}
\newcommand{\wl}{w_{0,\lambda}^j}
\newcommand{\cl}{\gamma_{0,\lambda}^j}

\newcommand{\whs}{\wh{\sigma}}

\newcommand{\dbarado}{\dbar^*_{h_0,\omega_0}}

\newcommand{\Sf}{S^q_{f,\omega,h}}

\newcommand{\into}{\int_{X_0}}

\newcommand{\pz}{\partial z}
\newcommand{\delho}{\partial_{h_0}}
\newcommand{\delhy}{\partial_{h_y}}
\newcommand{\delg}{\partial_g}

\newcommand{\lara}[2]{\langle{#1},{#2}\rangle}

\newcommand{\iO}[1]{i\Theta_{#1}}

\newcommand{\gaus}[1]{\lfloor {#1}\rfloor}


\begin{document}
\title
[Positivity of higher direct images] 
{Positivity of higher direct images for pseudo-effective line bundles}
\author{Yuta Watanabe}
\address{Department of Mathematics, Faculty of Fundamental Science and Engineering, Chuo University.
1-13-27 Kasuga, Bunkyo-ku, Tokyo 112-8551, Japan}
\email{{\tt wyuta.math@gmail.com}, {\tt wyuta@math.chuo-u.ac.jp}}

\begin{abstract}
    This paper investigates higher direct image sheaves twisted by pseudo-effective line bundles with their multiplier ideal sheaves. 
    Over a Zariski open subset of the base, these higher direct image sheaves are locally free and admit suitable Hodge-type metrics that are smooth and Nakano semi-positive. 
    Furthermore, we establish that any extension of these metrics as singular Hermitian metrics on torsion-free coherent sheaves satisfies Nakano semi-positivity and the minimal extension property, 
    and that these metrics extend uniquely with respect to Griffiths semi-positivity.
\end{abstract}


\maketitle





\section{Introduction}

The study of positivity for direct images of adjoint bundles goes back to the work of Griffiths \cite{Gri70} on variations of Hodge structure.
Fujita~\cite{Fuj78}, Kawamata~\cite{Kaw81}, Viehweg~\cite{Vie83}, and Koll\'ar~\cite{Kol86a,Kol86b}, among others, subsequently generalized these results in more algebro-geometric settings. 
Moriwaki~\cite{Mor87} and Fujino~\cite{Fuj04} also obtained positivity results for higher direct images.

From the viewpoint of complex analysis, the study of positivity for direct image sheaves began with work on the variation of Bergman kernels by Maitani and Yamaguchi \cite{MY04}. 
Berndtsson \cite{Ber09} developed this approach further and showed that the direct image sheaf inherits Nakano positivity.
To prove this, Berndtsson computed the curvature of the Hodge metric defined by fiber integration. 
The resulting formula shows that both the curvature of the original bundle and the Kodaira--Spencer classes describing the variation of the fibers' complex structures contribute to positivity.
Mourougane and Takayama \cite{MT08,MT09} subsequently generalized this method to show that higher direct image sheaves also inherit Nakano semi-positivity.
These analytic results concern smooth Hermitian metrics.

Turning to singular Hermitian metrics, Demailly \cite{Dem93} showed how such metrics can be used to study algebro-geometric notions of positivity for a line bundle $L$ by analytic means. 
Over the past decade, the theory of singular Hermitian metrics on vector bundles and, more generally, on torsion-free sheaves has developed considerably in complex geometry, e.g. \cite{BP08,PT18,HPS18,DNWZ23,Ina22}.
It is therefore natural to ask whether direct image sheaves inherit positivity in the sense of singular Hermitian metrics. 
More precisely, let $f:X\longrightarrow Y$ be a proper K\"ahler fiber space, and $L$ be a holomorphic line bundle on $X$ with a singular Hermitian metric $h$ satisfying a curvature positivity condition.
The question is whether the direct image sheaf $f_*(K_{X/Y}\otimes\lih)$ inherits a corresponding positivity property, where $\scr{I}(h)$ denotes the multiplier ideal sheaf associated with $h$.
Using the optimal form \cite{Blo13,GZ15} of the Ohsawa--Takegoshi $L^2$-extension theorem \cite{OT87}, the works \cite{PT18,HPS18} establish Griffiths semi-positivity of the natural singular Hermitian metric, 
i.e., the Hodge metric for $q=0$, on $f_*(K_{X/Y}\otimes\lih)$ when $(L,h)$ is pseudo-effective. In \cite{HPS18}, this positivity is derived using the minimal extension property.
More recently, the characterization of Nakano positivity in terms of optimal $L^2$-estimates of H\"ormander type provided the basis for a definition of singular Nakano positivity \cite{DNWZ23,Ina22}.
These optimal $L^2$-estimates have been used to establish singular Nakano positivity for direct image sheaves in the nef and big setting and in the pseudo-effective setting in \cite{WZ25} and \cite{IMW24}, respectively.
For related work, see \cite{BP08,CGP22,Tak24,Wat25}.

Earlier joint work of the author with Inayama and Matsumura \cite{IMW24} established Nakano positivity for direct image sheaves $f_*(K_{X/Y}\otimes\lih)$ in essentially full generality. 
We therefore aim to establish Nakano semi-positivity for the higher direct image sheaves $R^qf_*(K_{X/Y}\otimes\lih)$, which are torsion-free by \cite{Mat22}.
A natural approach would be to prove singular Nakano semi-positivity directly using optimal $L^2$-estimates, as in \cite{IMW24}. 
However, Professor Matsumura pointed out difficulties with this approach in the higher direct image setting. 
The author then revisited the earlier joint work with Zou \cite{WZ25} and the paper that inspired it \cite{CGP22}. 
This suggested that the curvature of the Hodge metric could be studied analytically when the singularities of $h$ are \textit{horizontal} with respect to $f$. 
The overall strategy is similar to that of \cite{Ber09,MT08}. 
We construct a suitable Hodge metric and derive a curvature formula optimized, from which we deduce the following result.

\begin{mainthm}\mainthmlabel{Main Theorem}
  Let $f:X\longrightarrow Y$ be a holomorphic map of complex manifolds, which is smooth, proper, K\"{a}hler, surjective, and with connected fibers. 
  Let $L$ be a holomorphic line bundle on $X$ with a singular Hermite metric $h$ such that $\iO{L,h}\geq0$ in the sense of currents. Assume that $h$ has analytic singularities modulo smooth functions. 
  Let $q\geq0$ and $Y^\circ_{q,h}$ denote the maximal Zariski open subset of $Y$ on which the higher direct image sheaf $\hdi$ is locally free. 
  Then, there exists a Zariski open subset $Y^\circ\subset Y^\circ_{q,h}$ satisfying the following properties: 
  \begin{itemize}
    \item The holomorphic vector bundle $\hdi|_{Y^\circ}$ admits a smooth Hodge metric $g$ that is Nakano semi-positive. 
    \item Any singular Hermitian metric $\breve{g}$ on $\hdi$ obtained by extending $g$ with arbitrary measurable Hermitian forms on $Y^\circ_{q,h}\setminus Y^\circ$ is Nakano semi-positive and has the minimal extension property.
    \item The smooth Hodge metric $g$ over $Y^\circ$ extends uniquely to a Griffiths semi-positive singular Hermitian metric $\tl{g}$ on $\hdi$.
  \end{itemize}
\end{mainthm}

This theorem answers \cite[Question~1.3-2]{Tak24} and establishes a stronger positivity property, namely Nakano semi-positivity in the appropriate sense. 
Even when $q=0$, it gives a more precise description of positivity than that obtained in \cite{IMW24}, though it requires that $h$ have analytic singularities.
Together with the underlying curvature formula (see Theorem~\ref{Theorem: curvature formula}), it generalizes the corresponding results of \cite{MT08,CGP22}.

Most of this paper is devoted to establishing the first assertion.
As explained in $\S$\ref{subsection: Occurrence horizontal setting}, the Zariski open subset $Y^\circ$ is chosen such that, after taking a log resolution $\pi:\widetilde X\longrightarrow X$ of the singularities of $h$, 
the resulting simple normal crossing divisor has relative simple normal crossings with respect to the Kähler morphism $\widetilde f:=f\circ\pi$ over $Y^\circ$. 
In other words, we mainly work in the setting where, after taking a log resolution, the multiplier ideal sheaf $\mathscr I(\pi^*h)$ becomes invertible 
and the singularities of $\pi^*h$ are \textit{horizontal} with respect to the Kähler morphism $\widetilde f$, namely, the horizontal singularity setting (\ref{Setting: horizontal singularity}).
In this setting, $\widetilde f$ is smooth over $Y^\circ$, and $\mathscr I(\pi^*h)|_{\widetilde X_y}=\mathscr I(\pi^*h|_{\widetilde X_y})$ for every $y\in Y^\circ$, where $\widetilde X_y:=\widetilde f^{-1}(y)$.
Furthermore, in local coordinates adapted to the relative simple normal crossings divisor, the singular part of the weight of $\pi^*h$ is independent of the base coordinates. 
This allows us to show that the Hodge metric constructed in this setting is \textit{smooth} on $Y^\circ$ (see Remark \ref{Remark: smoothness of Hodge metric g}).
Our choice to work in this setting (\ref{Setting: horizontal singularity}) was inspired by \cite{CGP22,WZ25}.
Here, Grauert--Riemenschneider-type vanishing \cite[Corollary~1.5]{Mat22} and the projection formula give $R^q\pi_*(K_{\tl{X}/X}\otimes\pi^*L\otimes\scr{I}(\pi^*h))=0$ for $q>0$.
Therefore, for every $q\geq0$, the Leray spectral sequence then gives
\[
  R^qf_*(K_{X/Y}\otimes\lih)\cong R^q\tl{f}_*(K_{\tl{X}/Y}\otimes\pi^*L\otimes\scr{I}(\pi^*h)).
\]

Our approach is to extend the Mourougane--Takayama framework \cite{MT08}, which relies on the harmonic theory developed by Takegoshi \cite{Tak95,Tak97}, to the pseudo-effective case.
Building on work on the hard Lefschetz theorem for pseudo-effective line bundles \cite{DPS01,Wu21}, the recent work \cite{QZ26} extends Takegoshi's harmonic theory to this setting. 
We first use this extension to construct the Hodge metric on the higher direct image sheaf. 
Guided by the curvature formula in \cite{MT08,CGP22}, we then establish a refined curvature formula from which Nakano semi-positivity follows.


\section{Preliminaries}

\subsection{Set-up}\label{subsec:setup}

We introduce the geometric setting used throughout this paper, beginning with the assumptions common to (\ref{Setting: local}) and (\ref{Setting: horizontal singularity}).

Let $f:X\longrightarrow Y$ be a smooth, proper, surjective K\"ahler morphism with connected fibers, where $X$ and $Y$ are complex manifolds of dimensions $n+m$ and $m$, respectively.
Let $L$ be a holomorphic line bundle on $X$ equipped with a singular Hermitian metric $h$ such that $\iO{L,h}\geq0$ in the sense of currents. 
In this case, $(L,h)$ is said to be \textit{singular semi}-\textit{positive}. In particular, if $X$ is compact, then $L$ is pseudo-effective.

Let $\omega_f$ be a relative K\"ahler metric for $f$, and denote its de Rham cohomology class by $\kappa_f:=\{\omega_f\}$.
For each $y\in Y$, we write $X_y:=f^{-1}(y), \omega_y:=\omega_f|_{X_y}, (L_y,h_y):=(L,h)|_{X_y}$.
Then $(X_y,\omega_y)$ is a compact K\"ahler manifold, and $(L_y,h_y)$ is pseudo-effective whenever the restriction $h_y$ is well-defined as a singular Hermitian metric.

\medskip
\refstepcounter{setting}
\noindent
\textbf{(I) The local setting.}\label{Setting: local}
We consider the above setting over a sufficiently small coordinate ball $Y\subset\mathbb{C}^m$, with holomorphic coordinates $t=(t_1,\ldots,t_m)$.
After shrinking $Y$ if necessary, we may choose a constant $c>0$ such that $\omega=\omega_f+cf^*(i\sum_jdt_j\wedge d\overline{t}_j)$ is a K\"ahler form on $X$.
This modification preserves the relative K\"ahler class $\kappa_f$ and satisfies $\omega|_{X_y}=\omega_y$ for every $y\in Y$.

We fix a system of local coordinates adapted to $f$.
Since $f$ is a holomorphic submersion, after restricting to a sufficiently small coordinate neighborhood $W\subset Y$, there exists a locally finite covering $\mathcal{U}=\{U_\alpha\}_\alpha$ of $f^{-1}(W)$
such that each $U_\alpha$ admits holomorphic coordinates $(z_\alpha,t)=(z_\alpha^1,\ldots,z_\alpha^n,t_1,\ldots,t_m)$ identifying $U_\alpha$ with a product $D_\alpha\times W$, where $D_\alpha\subset\mathbb{C}^n$ is a domain, and $f|_{U_\alpha}$ corresponds to the projection onto $W$.

On every nonempty intersection $U_\alpha\cap U_\beta$,
the coordinate changes take the form $z_\alpha^j=F_{\alpha\beta}^{j}(z_\beta,t)$ for $1\leq j\leq n$ where each $F_{\alpha\beta}^{j}$ is holomorphic in both $z_\beta$ and $t$.
We refer to such a covering $\mathcal{U}$ as a system of admissible charts over $W$.

\medskip
\refstepcounter{setting}
\noindent
\textbf{(II) The horizontal singularity setting.}\label{Setting: horizontal singularity}
We assume that the singularities of $h$ are supported on a divisor $D=\sum_{\ell=1}^{N}D_\ell$ on $X$ satisfying the following conditions.

\begin{itemize}
\item The divisor $D$ has relative simple normal crossings over $Y$. More precisely, for every point $x\in X$, there exist local coordinates $(z,t)=(z_1,\ldots,z_n,t_1,\ldots,t_m)$ centered at $x$ and adapted to $f$ such that $f(z,t)=t$ and $D=\{z_1\cdots z_r=0\}$ locally, for some $0\leq r\leq n$.
\item The singular Hermitian metric $h$ has divisorial singularities along $D$. Namely, in the above coordinates, its local weight $\varphi$ is given by
\[
  \varphi(z,t)=\sum_{j=1}^{r}a_j\log|z_j|^2+\psi(z,t),
\]
where $a_j\geq 0$ are constants and $\psi$ is a smooth real-valued function.
\end{itemize}

In particular, the singularities of $h$ are horizontal with respect to $f$, and the restriction $h_y=h|_{X_y}$ is well-defined for every $y\in Y$.
Moreover, the multiplier ideal sheaves satisfy
\begin{align*}
  \scr{I}(h)|_{X_y}=\scr{I}(h_y).
\end{align*}

Since most of our arguments are local on the base, we shall mainly work in setting (\ref{Setting: local}), together with setting (\ref{Setting: horizontal singularity}) whenever the horizontal singularity assumptions are required.

\subsection{Occurrence of the horizontal singularity setting}\label{subsection: Occurrence horizontal setting}
The horizontal singularity setting (\ref{Setting: horizontal singularity}) arises naturally through log resolutions, at least over a dense Zariski open subset of the base, provided that $h$ has analytic singularities modulo smooth functions.
Here, we say that $h$ has \textit{analytic singularities modulo smooth functions} if its local weights can be written as
\[
  \varphi=c\log\Big(\sum_{j=1}^{N}|g_j|^2\Big)+\psi,
\]
where $c>0$ is a real constant, the $g_j$ are holomorphic functions, and $\psi$ is smooth.

Indeed, let $\pi:\tx\longrightarrow X$ be a log resolution of the singularities of $h$.
Then the pullback metric $\pi^*h$ has divisorial singularities modulo smooth functions along a simple normal crossing divisor $\widetilde D$.
Furthermore, after restricting the base to a suitable dense Zariski open subset $Y^\circ\subset Y$, the induced morphism
\[
  \tl{f}:=f\circ\pi:\tl{f}^{-1}(Y^\circ)\longrightarrow Y^\circ
\]
is smooth, and $\widetilde D$ has relative simple normal crossings over $Y^\circ$. 

Hence, the horizontal singularity setting (\ref{Setting: horizontal singularity}) is satisfied for the resolved family over $Y^\circ$.
We refer to \cite[\S 2]{CGP22}, \cite{WZ25} for a related discussion.

\subsection{Notation}

For a holomorphic line bundle $L$ with a singular Hermitian metric $h$, we denote by $D_h=\partial_h+\bar\partial$ the Chern connection on the locus where $h$ is smooth and by $\Theta_{L,h}$ its curvature current.
Locally, writing $h=e^{-\varphi}$ in a holomorphic frame, we have
\[
  \partial_h=\partial-\partial\varphi\wedge\,\,, \qquad \iO{L,h}=\idd\varphi.
\]

For an open subset $U\subset X$, we denote by $A^{p,q}(U,L)$ the space of smooth $L$-valued $(p,q)$-forms on $U$, and by $L^2_{p,q}(U,L,h,\omega)$ (resp.\ $L^{2,loc}_{p,q}(U,L,h,\omega)$)
the space of square-integrable (resp.\ locally square-integrable) $L$-valued $(p,q)$-forms with respect to $h$ and $\omega$.
We omit $\omega$ from the notation whenever it is clear from the context or the underlying manifold is compact, in which case the $L^2$-space is independent of the choice of $\omega$.

Denote by $c_d:=i^{d^2}=i^d(-1)^{d(d-1)/2}$ for any non-negative integer $d$. 
Let $t=(t_1,\ldots,t_m)$ be the coordinates on $\bb{C}^m$.  
\begin{itemize}
  \item $dt:=dt_1\wedge\cdots\wedge dt_m$, $d\overline{t}:=\overline{dt}$, and $dV_t:=c_m\,dt\wedge d\bar t=\bigwedge_{j=1}^midt_j\wedge d\overline{t_j}>0$. 
  \item Let $\widehat{dt_j}$ be a $(m-1,0)$-form without $dt_j$ such that $dt_j\wedge\widehat{dt_j}=dt$, and $\widehat{d\overline{t_j}}=\overline{\widehat{dt_j}}$. 
  \item Let $\widehat{dt_j\wedge d\overline{t_k}}$ be a $(m-1,m-1)$-form without $dt_j$ and $d\overline{t_k}$ such that $idt_j\wedge d\overline{t_k}\wedge \widehat{dt_j\wedge d\bar t_k}=c_mdt\wedge d\overline{t}$. 
\end{itemize}
We use the same notation for their pullbacks to $X$, often omitting $f^*$, as in writing $dt$ instead of $f^*dt$.

\subsection{Relative differential forms}

Let $f:(X,\omega)\longrightarrow Y\subset\bb{C}^m$ be as in local setting (\ref{Setting: local}). 
For a subset $S\subset Y$, we denote by $X_S=f^{-1}(S)$ and $L_S=L|_S$. We briefly recall the notion of relative differential forms following \cite[\S 3.1]{MT08}.

Let $U\subset X$ be an open subset. For $u,v\in A^{p,q}(U,L)$, we write $u\sim v$ if $u|_{X_y\cap\,U}=v|_{X_y\cap\,U}$ for every $y\in Y$.
We define the space of relative differential forms by 
\[
A^{p,q}(U/Y,L):=A^{p,q}(U,L)/\sim,
\] 
and denote the class represented by $u$ by $[u]$.
In local coordinates $(z,t)$ adapted to $f$, every $u\in A^{p,q}(U,L)$ can be written as $u=\sum u_{IJ}\,dz^I\wedge d\bar z^J+R$,
where each term of $R$ contains at least one $dt_j$ or $d\bar t_j$. The class $[u]$ is thus determined by the coefficients $u_{IJ}$.

A class $[u]\in A^{p,0}(U/Y,L)$ is said to be \textit{holomorphic on each fiber} (resp. \textit{relatively holomorphic}) if the coefficients $u_I$ are holomorphic in $z$ (resp. in $(z,t)$).
Equivalently, $[u]$ is relatively holomorphic if and only if $u\wedge f^*dt$ is locally holomorphic over each coordinate chart of $Y$. 
For an open subset $W\subset Y$, we set
\[
\begin{aligned}
  A^0(W,f_*\Omega^p_{X/Y}(L))
  &:=\left\{[u]\in A^{p,0}(X_W/W,L)\ \middle|\ \text{$[u]$ is holomorphic on each fiber}\right\},\\
  H^0(W,f_*\Omega^p_{X/Y}(L))
  &:=\left\{[u]\in A^{p,0}(X_W/W,L)\ \middle|\ \text{$[u]$ is relatively holomorphic}\right\}.
\end{aligned}
\]
They are modules over $A^0(W,\mathbb C)$ and $\mathcal O_Y(W)$, respectively, with $\alpha[u]:=[(f^*\alpha)u]$.

Under the Setting~(\ref{Setting: horizontal singularity}),
we have $\mathscr I(h)|_{X_y}=\mathscr I(h_y)$ for every $y\in Y$.
We define 
\begin{align*}
  A^0(W,f_*(\Omega^p_{X/Y}\!\otimes\! L\!\otimes\!\mathscr I(h)))
  &\!:=\!\left\{[u]\!\in\! A^0(W,f_*\Omega^p_{X/Y}(L))\middle|\,u|_{X_y}\!\!\in\! L^2_{p,0}(X_y,L_y,h_y),\! \forall y\!\in\! W\right\}\!,\\ 
  H^0(W,f_*(\Omega^p_{X/Y}\!\otimes\! L\!\otimes\!\mathscr I(h)))
  &\!:=\!\left\{[u]\!\in\! H^0(W,f_*\Omega^p_{X/Y}(L))\middle|\,u|_{X_y}\!\!\in\! L^2_{p,0}(X_y,L_y,h_y),\! \forall y\!\in\! W\right\}\!.
\end{align*}
Indeed, the fiberwise $L^2$-integrability implies that each holomorphic coefficient $u_I$ is divisible by a local generator of $\mathscr I(h)$ in coordinates adapted to $f$ and $D$.

Equivalently, a class $[u]\in A^{p,0}(X_W/W,L)$ belongs to $H^0(W,f_*(\Omega^p_{X/Y}\otimes L\otimes\mathscr I(h)))$ if and only if $u\wedge f^*dt$ is holomorphic and locally $L^2$-integrable with respect to $h$ over each coordinate chart of $W$.

\begin{remark}[{\cite[Remark 3.2]{MT08}}]\label{Remark: MT08, 3.2}
  For $u,v\in A^{p,q}(X_W,L)$, a relation $u\wedge dt=v\wedge dt$ implies $[u]=[v]$ in $A^{p,q}(X_W/W,L)$, and the converse holds true in case $q = 0$.
\end{remark}

\section{Harmonic theory for pseudo-effective line bundles}

\subsection{Absolute setting}\label{subsection: Absolute setting}

Let $(X_0,\omega_0)$ be an $n$-dimensional compact Kähler manifold, and $(L_0,h_0)$ be a pseudo-effective line bundle on $X_0$.
By the Hard Lefschetz theorem for pseudo-effective line bundles \cite[Theorem 2.1.1]{DPS01}, the Lefschetz operator 
\[
L_0^q:H^0(X_0,\Omega_{X_0}^{n-q}\otimes L_0\otimes\scr{I}(h_0))\longrightarrow H^q(X_0,K_{X_0}\otimes L_0\otimes\scr{I}(h_0))
\]
induced by $\omega_0^q\wedge\bullet$ is surjective. Furthermore, by \cite[Theorem 2(a)]{Wu21}, its restriction
\[
L_0^q:H^0(X_0,\Omega_{X_0}^{n-q}\otimes L_0\otimes\scr{I}(h_0))\cap\Ker\partial_{h_0}\overset{\simeq}{\longrightarrow}H^q(X_0,K_{X_0}\otimes L_0\otimes\scr{I}(h_0))
\]
is an isomorphism, where $\partial_{h_0}$ is understood in the sense of currents.

We define the space of $L^2$-harmonic $(n,q)$-forms by
\[
\begin{aligned}
\mathcal H^{n,q}(X_0,L_0,h_0):=\bigl\{u\in L^2_{n,q}(X_0,L_0,h_0)\mid\dbar u=0,\, \dbar^*_{h_0,\omega_0}u=0\bigr\},
\end{aligned}
\]
where $\dbar^*_{h_0,\omega_0}$ denotes the Hilbert adjoint of $\dbar_{max}$.
On the other hand, by \cite[Theorems 1.1 and 3.3]{QZ26}, the natural map induces an isomorphism
\[
\iota_0:\mathcal H^{n,q}(X_0,L_0,h_0) \overset{\simeq}{\longrightarrow}H^q(X_0,K_{X_0}\otimes L_0 \otimes\mathscr I(h_0)),
\]
and the Hodge $*$-operator yields an isomorphism
\[
*_0:\mathcal H^{n,q}(X_0,L_0,h_0)\overset{\simeq}{\longrightarrow}H^0(X_0,\Omega_{X_0}^{n-q}\otimes L_0\otimes\mathscr I(h_0))\cap\Ker\partial_{h_0}.
\]
Furthermore, the space of $L^2$-harmonic forms admits the following characterization: 
\begin{align*}
  \mathcal H^{n,q}(X_0,L_0,h_0)=\bigl\{u\in A^{n,q}(X_0,L_0)\cap L^2_{n,q}(X_0,L_0,h_0)\mid \dbar u=0,\, \delho\!\ast_0u=0\bigr\}.
\end{align*}
Recalling that $\displaystyle\frac{c_{n-q}}{q!}\omega_0^q\wedge *_0u=u$ for $u\in\mathcal H^{n,q}(X_0,L_0,h_0)$, we have $\displaystyle\frac{c_{n-q}}{q!}L_0^q\circ *_0\circ\iota_0^{-1}=\mathrm{id}$.
Hence, we obtain the following direct sum decomposition:
\[
\begin{aligned}
H^0(X_0,\Omega_{X_0}^{n-q}\otimes L_0\otimes\mathscr I(h_0))&=\Ker L_0^q\oplus *_0\mathcal H^{n,q}(X_0,L_0,h_0)\\
&=\Ker L_0^q\oplus \big(H^0(X_0,\Omega_{X_0}^{n-q}\otimes L_0\otimes\scr{I}(h_0))\cap\Ker\partial_{h_0}\big).
\end{aligned}
\]

These results extend the harmonic theory of Enoki \cite{Eno93} and Takegoshi \cite{Tak95,Tak97} to pseudo-effective line bundles endowed with singular Hermitian metrics.

\subsection{Localized relative setting} 

Let $f:(X,\omega)\longrightarrow Y\subset\bb{C}^m$ be as in local setting (\ref{Setting: local}). 
We take a smooth plurisubharmonic exhaustion function $\Phi=f^*\sum^m_{j=1}|t_j|^2$ on $X$, where $X$ is holomorphically convex. 
Following \cite{QZ26}, we set the following subspace of $L$-valued harmonic $(n+m,q)$-forms with respect to $\omega$ and $h$:
\begin{align*}
    \cal{H}^{n+m,q}(X,L,h,\Phi)\!:=\Bigl\{u\!\in\! A^{n+m,q}\cap L^{2,loc}_{n+m,q}(X,\!L,h)\Big|\,\dbar u=0, \delh\!\ast\! u=0, \partial\Phi\!\wedge\!\ast u=0 \Bigr\},
\end{align*}

\begin{theorem}[{cf. \cite[Theorem 1.1, Theorem 3.3 and Corollary 3.7]{QZ26}}]\label{Theorem: QZ26}
    
    \noindent
    \begin{itemize}
      \item [(1)] The space $\cal{H}^{n+m,q}(X,L,h,\Phi)$ does not depend on smooth plurisubharmonic exhaustion functions $\Phi$, and represents $H^q(X,K_X\otimes L\otimes\scr{I}(h))$ as a torsion-free $\cal{O}(Y)$-module. 
      \item [(2)] We have the following commutative diagram: 
    \[
    \begin{tikzcd}[
  row sep=-2pt,
  column sep=normal
]
    \cal{H}^{n+m,q}(X,L,h,\Phi) \arrow[drrr, "\simeq", "\ast"']  \arrow[dd, "\simeq", "\iota"']  &&&  \\  
     &&&
    \begin{gathered}
        H^0(X,\Omega_X^{n+m-q}\otimes L\otimes\scr{I}(h)) \\ 
        \hspace{6mm}\cap\,\Ker\delh\cap\Ker(\partial\Phi\wedge\bullet)
    \end{gathered} \arrow[dlll, "\frac{c_{n+m-q}}{q!}\omega^q\wedge\bullet", "\simeq"'] \\
    H^q(X,K_X\otimes L\otimes\scr{I}(h))
    \end{tikzcd}
    \] 
    In particular, the Hodge $\ast$-operator and the Lefschetz homomorphism 
    \[
    \begin{tikzcd}[row sep=0pt, column sep=normal]
      \hspace{12mm} L^q_{\omega,h}:=\omega^q\wedge\bullet:H^0(X,\Omega_X^{n+m-q}\otimes L\otimes\scr{I}(h)) \arrow[r,two heads] & H^q(X,K_X\otimes L\otimes\scr{I}(h))
    \end{tikzcd}
    \] 
    satisfy $(c_{n+m-q}/q!)L^q_{\omega,h}\circ\ast\circ\iota^{-1}=\rom{id}$.
      \item [(3)] Let $u\in\calH$. Then the form $\ast u\in H^0(X,\Omega_X^{n+m-q}\otimes L\otimes\scr{I}(h))$ is saturated in base variables, i.e., $\ast u=\sigma_u\wedge dt$ for some $[\sigma_u]\in H^0(X,\Omega_{X/Y}^{n-q}\otimes L\otimes\scr{I}(h))$. 
      In particular, $u=(c_{n+m-q}/q!)\omega^q\wedge\sigma_u\wedge dt$ and the map $u\mapsto[\sigma_u]$ is well-defined. Thus the Hodge $\ast$-operator induces an injective homomorphism 
    \[
    \begin{tikzcd}[row sep=0pt, column sep=normal]
      \hspace{8mm} S^q_{\omega,h}:\calH \arrow[r,hook] & H^0(X,\Omega_{X/Y}^{n-q}\otimes L\otimes\scr{I}(h)).
    \end{tikzcd}
    \] 
    \end{itemize}
\end{theorem}

\begin{proof}
(1) and (2) follow from \cite[Theorems 1.1 and 3.3]{QZ26}.

(3) By \cite[Corollary 3.7 (2)]{QZ26}, we have $\ast u=\sigma_u\wedge dt$ for some $[\sigma_u]\in H^0(X,\Omega_{X/Y}^{n-q}\otimes L)$.
Since the holomorphic coefficients of $\ast u$ belong to $\scr{I}(h)$, the same holds for $\sigma_u$. Hence, $[\sigma_u]\in H^0(X,\Omega_{X/Y}^{n-q}\otimes L\otimes\scr{I}(h))$.
To show that the map is well-defined, suppose that $\sigma_u,\tl{\sigma}_u\in A^{p,0}(X,L)$ satisfy $\ast u=\sigma_u\wedge dt=\tl{\sigma}_u\wedge dt$.
Then $(\sigma_u-\tl{\sigma}_u)\wedge dt=0$. By the definition of relative differential forms, we obtain $[\sigma_u-\tl{\sigma}_u]=0$, that is, $[\sigma_u]=[\tl{\sigma}_u]$. Thus, the map $u\mapsto[\sigma_u]$ is well-defined.
Finally, we prove injectivity. Suppose that $S^q_h(u)=0$. Then $[\sigma_u]=0$, and hence $\sigma_u=\sum_{j=1}^m\eta_j\wedge dt_j$.
Consequently, $\ast u=\sigma_u\wedge dt=0$, which implies $u=\pm\ast\ast u=0$.
\end{proof}

Here, $\iota$ is induced by the natural quotient map $Z_{\dbar}^{n,q}(X,L,h)\longrightarrow H^q(X,K_X\otimes L\otimes\mathscr I(h))$ arising from the $L^2$-Dolbeault isomorphism.
In Theorem \ref{Theorem: QZ26} (3), we used the assumption that $f$ is smooth to decompose differential forms

We take a trivialization $\cal{O}_Y\overset{\simeq}{\longrightarrow}K_Y$ given by $1\mapsto dt$, which induces isomorphisms of sheaves $K_{X/Y}\cong K_{X/Y}\otimes f^*K_Y\cong K_X$ by $[u]\mapsto u\wedge dt$, and hence of cohomology groups 
$\alpha^q:H^q(X,K_{X/Y}\otimes L\otimes\scr{I}(h))\overset{\simeq}{\longrightarrow} H^q(X,K_X\otimes L\otimes\scr{I}(h))$.
We also have an injection $\Omega^{n-q}_{X/Y}\longrightarrow\Omega_X^{n+m-q}$ by $[\sigma]\mapsto\sigma\wedge dt$, and hence an injection $\beta^0: H^0(X,\Omega^{n-q}_{X/Y}\otimes L\otimes\scr{I}(h))\longrightarrow H^0(X,\Omega_X^{n+m-q}\otimes L\otimes\scr{I}(h))$. 
Combining with Theorem \ref{Theorem: QZ26}, we have 
\[
\begin{tikzcd}[row sep=0pt, column sep=large]
\iota^{-1}\circ\alpha^q:H^q(X,K_{X/Y}\otimes L\otimes\scr{I}(h))\arrow[r,"\simeq"] & \calH \\
\ast=\beta^0\circ S^q_{\omega,h}:\calH \arrow[r,hook] & H^0(X,\Omega_X^{n+m-q}\otimes L\otimes\scr{I}(h)) \\
(\alpha^q)^{-1}\circ L^q_{\omega,h}:H^0(X,\Omega_X^{n+m-q}\otimes L\otimes\scr{I}(h)) \arrow[r,two heads] & H^q(X,K_{X/Y}\otimes L\otimes\scr{I}(h)).
\end{tikzcd}
\]
Then Theorem \ref{Theorem: QZ26} (2) reads the following relative version: 

\begin{corollary}\label{Corollary: relative splitting}
  Let 
  \begin{align*}
    S^q_{f,\omega,h}&:=S^q_{\omega,h}\circ\iota^{-1}\circ\alpha^q:H^q(X,K_{X/Y}\otimes L\otimes\scr{I}(h))\longrightarrow H^0(X,\Omega_{X/Y}^{n-q}\otimes L\otimes\scr{I}(h)), \\
    L^q_{f,\omega,h}&:=(\alpha^q)^{-1}\circ L^q_{\omega,h}\circ\beta^0:H^0(X,\Omega_{X/Y}^{n-q}\otimes L\otimes\scr{I}(h))\longrightarrow H^q(X,K_{X/Y}\otimes L\otimes\scr{I}(h)).
  \end{align*}
  Then $(c_{n+m-q}/q!)L^q_{f,\omega,h}\circ S^q_{f,\omega,h}=\rom{id}$ on $H^q(X,K_{X/Y}\otimes L\otimes\scr{I}(h))$.
\end{corollary}

These homomorphisms are invariant under pullback perturbations of the Kähler form from the base.

\begin{lemma}\label{Lemma: independence under pullback}
Let $\omega_1$ and $\omega_2$ be \kah forms on $X$ such that $\omega_1-\omega_2=f^*\gamma$ for some smooth real $(1,1)$-form $\gamma$ on $Y$.
Then we have $L^q_{f,\omega_1,h}=L^q_{f,\omega_2,h}$ and $S^q_{f,\omega_1,h}=S^q_{f,\omega_2,h}$.
\end{lemma}

\begin{proof}
For any relative holomorphic form $[\sigma]\in H^0(X,\Omega^{n-q}_{X/Y}\otimes\lih)$, we have
\[
(\omega_1-\omega_2)\wedge\sigma\wedge dt=f^*\gamma\wedge\sigma\wedge dt=0.
\]
Consequently, $\omega_1^q\wedge\sigma\wedge dt=\omega_2^q\wedge\sigma\wedge dt$, which gives $L^q_{f,\omega_1,h}=L^q_{f,\omega_2,h}$.

Let $W:=\rom{Im}\,\ast_{\omega_j}=H^0(X,\Omega_X^{n+m-q}\otimes L\otimes\scr{I}(h))\cap\Ker\partial_h\cap\Ker(\partial\Phi\wedge\bullet)$, which is independent of $\omega_j$.
By Theorem \ref{Theorem: QZ26} (2), the maps $(c_{n+m-q}/q!)L^q_{\omega_j,h}|_W$  are isomorphisms. Moreover, by Theorem \ref{Theorem: QZ26} (3), every element of $W$ is of the form $\sigma\wedge dt$. 
Hence $L^q_{\omega_1,h}=L^q_{\omega_2,h}$ on $W$ by $L^q_{f,\omega_1,h}=L^q_{f,\omega_2,h}$, and thus $(L^q_{\omega_1,h})^{-1}=(L^q_{\omega_2,h})^{-1}$.
Since $\beta^0\circ S^q_{f,\omega_j,h}=(c_{n+m-q}/q!)^{-1}(L^q_{\omega_1,h})^{-1}\circ\alpha^q$, the injectivity of $\beta^0$ implies $S^q_{f,\omega_1,h}=S^q_{f,\omega_2,h}$.
\end{proof}

By \cite[Corollary 3.7 (4)]{QZ26}, these constructions can be carried out locally on \(Y\) and induce homomorphisms of direct image sheaves.

\begin{corollary}\label{Corollary: splitting sheaves}
  There exist homomorphisms induced from the Hodge $\ast$-operator and the Lefschetz homomorphism: 
  \begin{align*}
    S^q_{f,\omega,h}&:R^qf_*(K_{X/Y}\otimes L\otimes\scr{I}(h))\longrightarrow f_*(\Omega_{X/Y}^{n-q}\otimes\lih), \\
    L^q_{f,\omega,h}&:f_*(\Omega_{X/Y}^{n-q}\otimes\lih)\longrightarrow R^qf_*(K_{X/Y}\otimes\lih)
  \end{align*}
  such that $(c_{n+m-q}/q!)L^q_{f,\omega,h}\circ S^q_{f,\omega,h}=\rom{id}$ on $R^qf_*(K_{X/Y}\otimes\lih)$. In particular 
  \begin{align*}
    f_*(\Omega_{X/Y}^{n-q}\otimes\lih)=\cal{F}^{n-q}_{h,\omega}\oplus\cal{K}^{n-q}_{h,\omega}, \text{ with } \cal{F}^{n-q}_{h,\omega}=\rom{Im}\,S^q_{f,\omega,h} \text{ and } \cal{K}^{n-q}_{h,\omega}=\Ker L^q_{f,\omega,h}. 
  \end{align*}
\end{corollary}

\subsection{Local freeness}

Let $f:X\longrightarrow Y$ be a surjective proper Kähler morphism between complex manifolds, and let $(L,h)$ be a singular Hermitian line bundle with semipositive curvature currents.
By \cite[Corollary 1.5]{Mat22} and the projection formula, $R^qf_*(K_{X/Y}\otimes\lih)$ is torsion-free and coherent for every $q\geq 0$.
Furthermore, if $f$ is smooth and $h$ has only horizontal singularities as in (\ref{Setting: horizontal singularity}), then the following stronger result holds.

\begin{lemma}[{cf. \cite[Lemma 4.9]{MT08}}]\label{Lemma: horizontal base change}
Let $f:X\longrightarrow Y$ be a smooth, proper, and surjection Kähler morphism, and let $(L,h)$ be a pseudo-effective line bundle on $X$. 
Assume that $h$ satisfies the horizontal singularity setting (\ref{Setting: horizontal singularity}). 
Then, for every $q\geq 0$, the sheaf $R^qf_*(K_{X/Y}\otimes L\otimes\scr{I}(h))$ is locally free, and the natural base change map
\begin{align*}
  R^qf_*(K_{X/Y}\otimes L\otimes\scr{I}(h))\otimes_{\cal{O}_Y}\bb{C}(y)\overset{\simeq}{\longrightarrow} H^q(X_y,K_{X_y}\otimes L_y\otimes\scr{I}(h_y))
\end{align*}
is an isomorphism for every $y\in Y$.
\end{lemma}

\begin{proof}
  The proof proceeds as in \cite[Lemma 4.9]{MT08}, using Matsumura's torsion-freeness theorem.
Fix $y\in Y$. Since the assertion is local on $Y$, we may assume that $Y$ is a ball in $\bb{C}^m$ with coordinates $t=(t_1,\ldots,t_m)$ centered at $y$. Set 
\begin{align*}
  Y_j=\{t_1=\cdots=t_j=0\},\,\, X_j=f^{-1}(Y_j),\,\, f_j=f|_{X_j} \,\text{ and }\, \cal{G}_j=K_{X_j/Y_j}\otimes L_j\otimes\scr{I}(h_j),
\end{align*}
where $L_j=L|_{X_j}$ and $h_j=h|_{X_j}$.
By the horizontal singularity setting (\ref{Setting: horizontal singularity}),  $\scr{I}(h_j)=\scr{I}(h)|_{X_j}$ is invertible, and $\cal{G}_j|_{X_{j+1}}\simeq\cal{G}_{j+1}$.
Hence we have exact sequences 
\[
\begin{tikzcd}
0 \arrow[r]
& \cal{G}_j \arrow[r, "\cdot f_j^*t_{j+1}"]
& \cal{G}_j \arrow[r]
& \iota_{j*}\cal{G}_{j+1} \arrow[r]
& 0,
\end{tikzcd}
\]
where $\iota_j:X_{j+1}\hookrightarrow X_j$.
By Matsumura's torsion-freeness theorem \cite[Corollary 1.5]{Mat22}, the sheaf $R^qf_{j*}\cal{G}_j$ is torsion-free for every $q$.
Thus multiplication by $t_{j+1}$ on $R^{q+1}f_{j*}\cal{G}_j$ is injective. The associated long exact sequences therefore yield surjective restriction maps $R^qf_{j*}\cal{G}_j\longrightarrow R^qf_{j+1*}\cal{G}_{j+1}$.
Iterating these maps gives a surjection $(R^qf_*\cal{G}_0)_y\longrightarrow H^q(X_y,\cal{G}_m)$. Since $\cal{G}_0$ is locally free and $f$ is smooth, $\cal{G}_0$ is flat over $Y$.
The conclusion now follows from cohomology and base change \cite[Lemma 4.8]{MT08}.
\end{proof}

\section{Properties in the horizontal singularity setting}

In this section, we work under local setting (\ref{Setting: local}) and assume that $h$ satisfies the horizontal singularity setting (\ref{Setting: horizontal singularity}). 
We show that the results analogous to \cite[Lemmas 4.5, 4.7 and 3.1]{MT08} also hold in this setting.
Let $[\sigma]\in H^0(Y,f_*(\Omega^p_{X/Y}\otimes L))=H^0(X,\Omega^p_{X/Y}\otimes L)$.
Let $\{(U_\alpha;z_\alpha,t)\}_{\alpha\in\Lambda}$ be an admissible charts over $Y$. 
For each \(\alpha\), there exists a unique holomorphic representative $\sigma_\alpha\in H^0(U_\alpha,\Omega_X^p\otimes L)$ of \([\sigma]|_{U_\alpha}\) containing no \(dt_j\)-components. 
More precisely, $[\sigma_\alpha]=[\sigma]|_{U_\alpha}$ and $\sigma_\alpha=\sum_{|I|=p}\sigma_{\alpha,I}(z_\alpha,t)\,dz_I$.
Let \(\{\rho_\alpha\}\) be a smooth partition of unity subordinate to \(\{U_\alpha\}\). Then
\[
\sigma^\sharp:=\sum_\alpha\rho_\alpha\sigma_\alpha\in A^{p,0}(X,L)
\]
defines a globally smooth representative of \([\sigma]\), satisfying $[\sigma^\sharp]=[\sigma]$.
We call $\sigma^\sharp$ an \textit{admissible representative} of $[\sigma]$.

\begin{proposition}[{cf. \cite[Lemma 4.5]{MT08}}]\label{Proposition: key Prop. MT08, Lemma 4.5}
  Let $[\sigma]\in H^0(Y,\cal{F}^{n-q}_{h,\omega})$. 
  For an admissible representative \(\sigma^\sharp\in A^{n-q,0}(X,L)\cap L^{2,loc}_{n-q,0}(X,L,h)\) of \([\sigma]\) as constructed above, we obtain
  \begin{itemize}
    \item [(1)] $\delh\sigma^\sharp=\sum_j\mu^j\wedge dt_j$ in the sense of currents for some $\mu^j\!\in\! A^{n-q,0}\cap L^{2,loc}_{n-q,0}(X,L,h)$,
    \item [(2)] these $\mu^j$ are not unique, but $[\mu^j]\in A^{n-q,0}(X/Y,L)$ are well-defined for $\sigma^\sharp$, and
    \item [(3)] $\partial_{h_y}(\mu^j|_{X_y})=0$ in the sense of currents and $\mu^j|_{X_y}\in A^{n-q,0}\cap L^2_{n-q,0}(X_y,L_y,h_y)$. 
  \end{itemize}
\end{proposition}

\begin{proof}
  (1) There exists $\gamma\in H^0(Y,R^qf_*(K_{X/Y}\otimes\lih))=H^q(X,K_{X/Y}\otimes\lih)$ such that $[\sigma]=\Sf(\gamma)$ and $\alpha^q(\gamma)=\gamma\wedge dt\in H^q(X,K_X\otimes\lih)$.
  Set $u:=\iota^{-1}\circ\alpha^q(\gamma)=\iota^{-1}(\gamma\wedge dt)\in\calH$. 
  Then, by Theorem~\ref{Theorem: QZ26} (3), we have $[\sigma_u]=[\sigma]$ and $\ast u=\sigma_u\wedge dt=\sigma\wedge dt$.
  By Theorem~\ref{Theorem: QZ26} (2), we obtain $0=\delh\ast u=\delh(\sigma\wedge dt)=\delh\sigma\wedge dt$ in the sense of currents.
  Since $[\sigma^\sharp]=[\sigma]$, a similar argument shows that $\ast u=\sigma^\sharp\wedge dt$, and hence $\delh\sigma^\sharp\wedge dt=0$.
  Here, for each $\alpha$, there exists $\sigma_\alpha\in H^0(U_\alpha,\Omega_X^{n-q}\otimes\lih)$ such that $\ast u|_{U_\alpha}=\sigma_\alpha\wedge dt$ and $\sigma^\sharp=\sum_\alpha\rho_\alpha\sigma_\alpha$.

  By Setting (\ref{Setting: horizontal singularity}), on each $U_\alpha$, we have $D\cap U_\alpha=\{z_\alpha^1\cdots z_\alpha^{r_\alpha}=0\}$, and the local weight $\varphi_\alpha$ of $h$ can be written as
  $\varphi_\alpha=\sum^{r_\alpha}_{j=1}a_{\alpha,j}\log|z_\alpha^j|^2+\psi_\alpha$, where $a_{\alpha,j}\geq0$ and $\psi_\alpha\in C^{\infty}(U_\alpha)$.
  Here, the $(1,0)$-part of the Chern connection 
  \begin{align*}
    \delh|_{U_\alpha}=\partial-\partial\varphi_\alpha\wedge\, =\partial-\Big(\sum^{r_\alpha}_{j=1}a_{\alpha,j}\frac{dz_\alpha^j}{z^j_\alpha}+\partial\psi_\alpha\Big)\wedge\,
  \end{align*}
  is well-defined as an operator from smooth $L$-valued differential forms to currents with locally $L^1$ coefficients, and
  \begin{align*}
    \delh\sigma_\alpha=\partial_z\sigma_\alpha-\Big(\sum^{r_\alpha}_{j=1}a_{\alpha,j}\frac{dz_\alpha^j}{z^j_\alpha}+\partial_z\psi_\alpha\Big)\wedge\sigma_\alpha
    +\sum^m_{j=1}dt_j\wedge\Big(\partial_{t_j}\sigma_\alpha-\frac{\partial\psi_\alpha}{\partial t_j}\sigma_\alpha\Big).
  \end{align*}
  From $0=\delh\ast u|_{U_\alpha}=\delh\sigma_\alpha\wedge dt$, we first obtain $\partial_z\sigma_\alpha-\big(\sum^{r_\alpha}_{j=1}a_{\alpha,j}\frac{dz_\alpha^j}{z^j_\alpha}+\partial_z\psi_\alpha\big)\wedge\sigma_\alpha=0$ on $U_\alpha\setminus D$.
  Since each term is locally integrable, this equality also holds on $U_\alpha$ in the sense of currents.
  Thus, setting $\mu^j_\alpha:=(-1)^{n-q}(\partial_{t_j}\sigma_\alpha-(\psi_\alpha)_{t_j}\sigma_\alpha)$, we obtain $\delh\sigma_\alpha=\sum^m_{j=1}\mu^j_\alpha\wedge dt_j$.
  By $\sigma_\alpha\in H^0(U_\alpha,\Omega_X^{n-q}\otimes\lih)$, we have $\mu^j_\alpha\in A^{n-q,0}(U_\alpha,\bb{C})\cap L^{2,loc}_{n-q,0}(U_\alpha,h)$.

  Since $\delh$ is well-defined, the Leibniz rule yields $\delh\sigma^\sharp=\sum_\alpha(\partial\rho_\alpha\wedge\sigma_\alpha+\rho_\alpha\delh\sigma_\alpha)$ in the sense of currents. Moreover, for some constants $C,C'>0$, we obtain
  \begin{align*}
    |\partial\rho_\alpha\wedge\sigma_\alpha|^2_{h,\omega}\leq C|\partial\rho_\alpha|^2_\omega|\sigma_\alpha|^2_{h,\omega}\leq C'|\sigma_\alpha|^2_{h,\omega},
  \end{align*}
  Hence, $\partial\rho_\alpha\wedge\sigma_\alpha\in A^{n-q+1,0}\cap L^{2,loc}_{n-q+1,0}(U_\alpha,h)$, and consequently, $\delh\sigma^\sharp\in A^{n-q+1,0}\cap L^{2,loc}_{n-q+1,0}(X,L,h)$.
  Combining this with $\delh\sigma^\sharp\wedge dt=0$, we obtain (1).

  (2) By (1), we have $\mu^j\wedge dt=(-1)^{j-1}\delh\sigma^\sharp\wedge\wh{dt_j}$, which depends only on $\sigma^\sharp$.
  Thus, by Remark \ref{Remark: MT08, 3.2}, the class $[\mu^j]$ is well-defined for $\sigma^\sharp$.
  Indeed, if $\delh\sigma^\sharp=\sum_j\tl{\mu}^j\wedge dt_j$ is another such expression, then $0=(\mu^j-\tl{\mu}^j)\wedge dt$, and hence $[\mu^j]=[\tl{\mu}^j]\in A^{n-q,0}(X/Y,L)$. 

  (3) By the Setting (\ref{Setting: horizontal singularity}) and Fubini's theorem, we have $\mu^j|_{X_y}\in L^2_{n-q,0}(X_y,L_y,h_y)$.
  Since $h$ is smooth on $X\setminus D$, we have $\delh^2=0$ on $X\setminus D$. Hence, on $X\setminus D$, we obtain $0=\delh^2\sigma^\sharp=\delh\sum_j\mu^j\wedge dt_j=\sum_j\delh\mu^j\wedge dt_j$.
  Comparing the components containing exactly one $dt_j$, we obtain $\partial_{h_y}(\mu^j|_{X_y\setminus D_y})=0$, where we henceforth write $D_y:=D\cap X_y$.

  On the other hand, since $\mu^j$ is smooth, the expression
  \begin{align*}
    \partial_{h_y}(\mu^j|_{X_y})=\partial_z\mu^j-\Big(\sum^{r_\alpha}_{j=1}a_{\alpha,j}\frac{dz^j_\alpha}{z^j_\alpha}+\partial_z\psi_\alpha|_{X_y}\Big)\wedge\mu^j
  \end{align*}
  has locally $L^1_{loc}$ coefficients on each $U_\alpha$ with $U_\alpha\cap D_y\ne\emptyset$.
  Therefore, $\partial_{h_y}(\mu^j|_{X_y\setminus D_y})=0$ extends across $D\cap X_y$ in the sense of currents.
\end{proof}

Without the assumption that the singularities are horizontal, that is, of relative simple normal crossing type, 
it may not be possible to choose the forms $\mu^j$ to be smooth and $L^2$-integrable with respect to $h$. 
Since such a choice also requires local representatives $\sigma|_U$ with no $dt_j$-components, we introduced the admissible representative $\sigma^\sharp$.

Stalks or fibers at a point $y \in Y$ will be denoted by $f_*(\Omega^{n-q}_{X/Y}\otimes\lih)_y, (\cal{F}^{n-q}_{h,\omega})_y$ and $(\cal{K}^{n-q}_{h,\omega})_y$ respectively. 
Those stalks can be seen as subspaces of $H^0(X_y,\Omega^{n-q}_{X_y}\otimes L_y)$, i.e., $(\cal{F}^{n-q}_{h,\omega})_y\oplus(\cal{K}^{n-q}_{h,\omega})_y=f_*(\Omega^{n-q}_{X/Y}\otimes\lih)_y\subset H^0(X_y,\Omega^{n-q}_{X_y}\otimes L_y\otimes(\scr{I}(h)\cdot\cal{O}_{X_y}))$.

\begin{proposition}[{cf. \cite[Lemma 4.7]{MT08}}]\label{Proposition: key Prop. MT08, Lemma 4.7}
  Let $\sigma_y\in (\cal{F}^{n-q}_{h,\omega})_y$ and $\tau_y\in (\cal{K}^{n-q}_{h,\omega})_y$, and regard them as elements of $H^0(X_y,\Omega^{n-q}_{X_y}\otimes L_y\otimes\scr{I}(h_y))$. 
  Then
  \begin{itemize}
    \item [(1)] $\partial_{h_y}\sigma_y=0$ in the sense of currents on $X_y$,
    \item [(2)] $\omega^q_y\wedge\tau_y\in A^{n,q}(X_y,L_y)\cap L^2_{n,q}(X_y,L_y,h_y)$ is $\dbar$-exact, and
    \item [(3)] $(\sigma_y,\tau_y)_{h_y}=\int_{X_y}(c_{n-q}/q!)\omega^q_y\wedge\sigma_y\wedge h_y\overline{\tau_y}=0$. 
  \end{itemize}
\end{proposition}

\begin{proof}
  We will argue at $y=0$. (1) Since $Y$ is a unit ball in $\bb{C}^m$, there exists $[\sigma]\in H^0(Y,\cal{F}^{n-q}_{h,\omega})$ such that $\sigma|_{X_0}=\sigma_0$.
  By Proposition \ref{Proposition: key Prop. MT08, Lemma 4.5}, we choose a smooth representative $\sigma^\sharp\in A^{n-q,0}(X,L)$ of $[\sigma]$, which satisfies $[\sigma]=[\sigma^\sharp]$, 
  in particular $\sigma^\sharp|_{X_0}=\sigma|_{X_0}=\sigma_0$, and $\delh\sigma^\sharp=\sum_j\mu^j\wedge dt_j$ in the sense of currents for some $\mu^j\in A^{n-q,0}\cap L^{2,loc}_{n-q,0}(X,L,h)$.

  Thus, we obtain $\partial_{h_0}\sigma_0=0$. Indeed, since $\delh\sigma^\sharp|_{X_0\setminus D_0}=\sum_j\mu^j\wedge dt_j|_{X_0\setminus D_0}=0$ and $h$ is smooth on $X\setminus D$, we have $\partial_{h_0}\sigma_0|_{X_0\setminus D_0}=\delh\sigma^\sharp|_{X_0\setminus D_0}=0$.
  On the other hand, since $\sigma_0$ is smooth, as in Proposition \ref{Proposition: key Prop. MT08, Lemma 4.5} (3), the expression
  \begin{align*}
    \partial_{h_0}\sigma_0=\partial_z\sigma_0-\Big(\sum^{r_\alpha}_{j=1}a_{\alpha,j}\frac{dz^j_\alpha}{z^j_\alpha}+\partial_z\psi_\alpha|_{X_0}\Big)\wedge\sigma_0
  \end{align*}
  has locally $L^1$ coefficients on each $U_\alpha$ with $U_\alpha\cap D_0\ne\emptyset$.
  Therefore, the equality $\partial_{h_0}\sigma_0|_{X_0\setminus D_0}=0$ extends to $X_0$ in the sense of currents.

  (2) This follows by the same argument as in \cite[Lemma 4.7 (2)]{MT08}.
  In particular, the $L^2$-integrability with respect to $h$ also holds. Indeed, as is well known, $(\scr{L}^{n,\ast}_{L,h},\dbar)$ provides an $L^2$-Dolbeault resolution of $K_X\otimes\lih$, see \cite{Dem93}.
  Moreover, the Dolbeault complex with smooth coefficients $(\scr{A}^{n,\ast}_X(L)\otimes_{\cal{O}_X}\scr{I}(h),\dbar)$ also provides a resolution of $K_X\otimes\lih$. 
  Here, $\scr{A}_X^{n,q}$ denotes the sheaf of smooth $(n,q)$-forms. Thus, we may choose a smooth differential form as a representative. 

  (3) By (2), there exists $b\in A^{n,q-1}\cap L^{2,loc}_{n,q-1}(X,L,h)$ such that $\omega^q_0\wedge\tau_0=\dbar b_0$ on $X_0$, where $b_0=b|_{X_0}$.
  By the Setting (\ref{Setting: horizontal singularity}) and Fubini's theorem, we obtain $b_0\in A^{n,q-1}\cap L^2_{n,q-1}(X_0,L_0,h_0)$ and $\dbar b_0\in L^2_{n,q}(X_0,L_0,h_0)$.
  Moreover, since $\sigma_0\in H^0(X_0,\Omega^{n-q}_{X_0}\otimes L_0\otimes\scr{I}(h_0))$, the Cauchy--Schwarz inequality implies that $\sigma_0\wedge h_0\overline{\dbar b_0}$ is integrable on $X_0$, and hence we have $\int_{X_0}\omega^q_0\wedge\sigma_0\wedge h_0\overline{\tau_0}=\int_{X_0}\sigma_0\wedge h_0\overline{\dbar b_0}$.

  On $X_0\setminus D_0$, the smoothness of $h_0|_{X_0\setminus D_0}$ and (1) imply
  \begin{align*}
    \partial(\sigma_0\wedge h_0\overline{b_0})&= (\partial_{h_0}\sigma_0)\wedge h_0\overline{b_0}+(-1)^{n-q}\sigma_0\wedge h_0\overline{\dbar b_0}=(-1)^{n-q}\sigma_0\wedge h_0\overline{\dbar b_0},   \text{ and}\\
    \into\sigma_0\wedge h_0\overline{\dbar b_0}&=\int_{X_0\setminus D_0}\sigma_0\wedge h_0\overline{\dbar b_0}=(-1)^{n-q}\int_{X_0\setminus D_0}\partial(\sigma_0\wedge h_0\overline{b_0}).
  \end{align*}
  However, since $X_0\setminus D_0$ is open, we cannot directly apply Stokes' theorem. 
  To overcome this difficulty, we introduce a suitable cutoff function and show that the resulting boundary term tends to zero.

  We first examine the vanishing orders of $\sigma_0$ using its integrability. 
  On a local chart $U\subset X_0$ of $D_0=\{z_1\cdots z_r=0\}$, we write 
  $\varphi_0=\sum^r_{j=1}a_j\log|z_j|^2+\psi(z)$, where 
  $h_0|_U=e^{-\varphi_0}$, $a_j\geq0$ and $\psi$ is a smooth function.
  
  \begin{lemma}\label{Lemma: key vanishing orders}
    Locally, write $\varphi_0=\sum_{j=1}^r a_j\log|z_j|^2+\psi(z)$. Let $\sigma_0=\sum_{|I|=n-q}\sigma_I\,dz_I\in H^0(X_0,\Omega^{n-q}_{X_0}\otimes L_0\otimes\scr{I}(h_0))$ be a holomorphic section satisfying $\partial_{h_0}\sigma_0=0$.
    Then, for every $I$ and every $\ell\notin I$ with $a_\ell>\gaus{a_\ell}$, we have $z_\ell^{\gaus{a_\ell}+1}\mid\sigma_I$.
  \end{lemma}
  
\begin{proof}
  We can write $\sigma_I=\!\prod^r_{j=1}z_j^{\gaus{a_j}}s_I\!=\!z_\ell^{\gaus{a_\ell}}s_{I,\ell}$, where $s_I$ is also holomorphic and $s_{I,\ell}:=\prod^r_{\ell\ne j=1}z_j^{\gaus{a_j}}s_I$ for each $\ell\in\{1,\ldots,r\}$.
  By (1), we obtain
  \begin{align*}
    0=\partial_{h_0}\sigma_0&=\partial_z\sigma_0-\Big(\sum^r_{\ell=1}a_\ell\frac{dz_\ell}{z_\ell}+\partial_z\psi\Big)\wedge\sigma_0\\
    &=\sum^n_{j=1}\sum_{|I|=n-q}\Big(\frac{\partial\sigma_I}{\pz_j}-\frac{\partial\psi}{\pz_j}\sigma_I\Big)dz_j\wedge dz_I-\sum^r_{\ell=1}\sum_{|I|=n-q}\frac{a_\ell}{z_\ell}\sigma_Idz_\ell\wedge dz_I.
  \end{align*}
  Thus, for each $\ell\notin I$, comparing the coefficients of $dz_\ell\wedge dz_I$, we obtain
  \begin{align*}
    0=\frac{\partial\sigma_I}{\pz_\ell}-\frac{a_\ell}{z_\ell}\sigma_I-\frac{\partial\psi}{\pz_\ell}\sigma_I+\sum_{j\in I}\varepsilon_{I,j}\Big(\frac{\partial\sigma_{J_j}}{\pz_j}-\frac{\partial\psi}{\pz_j}\sigma_{J_j}\Big),
  \end{align*}
  where $J_j=(I\setminus\{j\})\cup\{\ell\}$ and $\varepsilon_{I,j}\in\{1,-1\}$. Consequently, we can write
  \begin{align*}
    0=(\gaus{a_\ell}-a_\ell)z_\ell^{\gaus{a_\ell}-1}s_{I,\ell}+z_\ell^{\gaus{a_\ell}}R_I,
  \end{align*}
  where $R_I$ is a smooth function. Hence, for any $\ell\notin I$, if $a_\ell>\gaus{a_\ell}$, then $z_\ell\mid s_{I,\ell}$, and thus $z_\ell^{\gaus{a_\ell}+1}\mid\sigma_I$.
\end{proof}

  Choose a cutoff function $\chi_\varepsilon=\prod^r_{j=1}\chi(|s_j|/\varepsilon)$ as in the following lemma \ref{Lemma: cutoff near D}. It vanishes near $D$, equals $1$ outside an $O(\varepsilon)$-neighborhood of $D$, and $\partial\chi_\varepsilon$ is locally a sum of terms $O(\varepsilon^{-1})dz_j$, each supported where $|z_j|\asymp\varepsilon$, up to bounded cutoff factors.
  By integrability and the Lebesgue dominated convergence theorem, we obtain
  \[
    \lim_{\varepsilon\to0}\int_{X_0}\chi_\varepsilon\sigma_0\wedge h_0\overline{\dbar b_0}=\into\sigma_0\wedge h_0\overline{\dbar b_0}.
  \]
  On the other hand, by Stokes' theorem, we have
  \[
    \into\chi_\varepsilon\sigma_0\wedge h_0\overline{\dbar b_0}=(-1)^{n-q}\into\chi_\varepsilon\partial(\sigma_0\wedge h_0\overline{b_0})=(-1)^{n-q+1}\into\partial\chi_\varepsilon\wedge\sigma_0\wedge h_0\overline{b_0}.
  \]
  Set $\delta_j=1$ if $a_j\!>\!\gaus{a_j}$ and $\delta_j=0$ if $a_j=\gaus{a_j}$. Recall that $b_0\in A^{n,q-1}\cap L^2_{n,q-1}(X_0,L_0,h_0)$.
  Locally, we have $|b_0|\leq C\prod^r_{j=1}|z_j|^{\gaus{a_j}}$. Combining this with the above divisibility $z_\ell^{\gaus{a_\ell}+1}\mid\sigma_I$ for any $\ell\notin I$ with $a_\ell>\gaus{a_\ell}$,
  we see that the local contribution to $\into\partial\chi_\varepsilon\wedge\sigma_0\wedge h_0\overline{b_0}$ is $O\big(\varepsilon^{1+\delta_j-2(a_j-\gaus{a_j})}\big)=o(1)$, since $1+\delta_j-2(a_j-\gaus{a_j})>0$.
  Integrability in the other variables and a finite cover therefore yield
  \begin{align*}
    \lim_{\varepsilon\to0}\into\partial\chi_\varepsilon\wedge\sigma_0\wedge h_0\overline{b_0}=0.
  \end{align*}

  In fact, by Lemma \ref{Lemma: key vanishing orders}, if $\ell\notin I$ then $|\sigma_I|\leq C|z_\ell|^{\gaus{a_\ell}+\delta_\ell}\prod_{j\ne \ell}|z_j|^{\gaus{a_j}}$. Locally, we have
  \begin{align*}
    ||O(\varepsilon^{-1})dz_\ell\wedge\sigma_0||^2_{L^2(h_0;S_{\ell,\varepsilon})}\leq C\varepsilon^{-2}\int_{r\varepsilon\leq|z_\ell|\leq R\varepsilon}|z_\ell|^{2(\delta_\ell-(a_\ell-\gaus{a_\ell}))}dV_{z_\ell}\leq C\varepsilon^{2(\delta_\ell-(a_\ell-\gaus{a_\ell}))}
  \end{align*}
  and $||O(1)\wedge\sigma_0||^2_{L^2(h_0;S_{\ell,\varepsilon})}\leq C\varepsilon^{2-2(a_\ell-\gaus{a_\ell})}$, where $e^{-\varphi_0}=e^{-\psi}\prod^r_{j=1}|z_j|^{-2a_j}$ and $S_{\ell,\varepsilon}:=\rom{supp}\,\partial\chi(|s_\ell|/\varepsilon)\subset\{r\varepsilon\leq|z_\ell|\leq R\varepsilon\}$.
  These yield $||\partial\chi_\varepsilon\wedge\sigma_0||^2_{L^2(h_0)}\!\leq\! C\sum^r_{j=1}\!\big(\varepsilon^{2(\delta_j-(a_j-\gaus{a_j}))}$ $+\varepsilon^{2-2(a_j-\gaus{a_j})}\big)\leq C$. 
  Hence, the Cauchy--Schwarz inequality implies 
  \begin{align*}
    \bigg|\into\partial\chi_\varepsilon\wedge\sigma_0\wedge h_0\overline{b_0}\bigg|&\leq C||\partial\chi_\varepsilon\wedge\sigma_0||_{L^2(h_0)}||b_0||_{L^2(\rom{supp}\,\partial\chi_\varepsilon,h_0)}\\
    &\leq C||b_0||_{L^2(\rom{supp}\,\partial\chi_\varepsilon,h_0)}\longrightarrow 0,
  \end{align*}
  if $\varepsilon\to0$.
  This proves $(\sigma_0,\tau_0)_{h_0}=(c_{n-q}/q!)\into\sigma_0\wedge h_0\overline{\dbar b_0}=0$.
\end{proof}

The proof of this proposition shows that, in the horizontal singularity setting (\ref{Setting: horizontal singularity}), Stokes' theorem applies as usual under the stated integrability conditions. 
This is a key observation. By developing this method further, we can carry out calculations analogous to those for smooth metrics and derive a corresponding curvature formula.

\begin{lemma}\label{Lemma: cutoff near D}
Let $s_j$ be the canonical section of $\mathcal O_X(D_j)$, and let $|s_j|$ denote its norm with respect to a smooth Hermitian metric.
Choose $\chi\in C^\infty([0,\infty),[0,1])$ such that $\chi=0$ on $[0,1]$ and $\chi=1$ on $[2,\infty)$, and set $\chi_\varepsilon=\prod_{j=1}^r\chi(|s_j|/\varepsilon)$.
Then the cutoff $\chi_\varepsilon$ vanishes near $D$ and converges pointwise to $1$ on $X\setminus D$.
Put $B_{j,\varepsilon}=\{\varepsilon\le |s_j|\le2\varepsilon\}$. In local coordinates where $D=\{z_1\cdots z_r=0\}$, we have 
\begin{align*}
  \partial\chi_\varepsilon=\sum_{j=1}^r\left(O(\varepsilon^{-1})\,dz_j+O(1)\right)\quad\text{and}\quad |\partial\chi_\varepsilon|\le\frac{C}{\varepsilon}\sum_{j=1}^r\mathbf 1_{B_{j,\varepsilon}},
\end{align*}
where each summand is supported in $B_{j,\varepsilon}$, up to bounded cutoff factors, and $B_{j,\varepsilon}\subset\{c\varepsilon\le |z_j|\le C\varepsilon\}$ for some positive constants $r,R$ on each coordinate neighborhood near $D$. 
Furthermore, suppose that $\varphi_0=\sum_{j=1}^r a_j\log|z_j|^2+\psi$, where $a_j>0$ and $\psi$ is smooth. 
Then, on every relatively compact adapted coordinate $U$,
\[
\begin{aligned}
\int_{B_{j,\varepsilon}\cap U}e^{\varphi_0}\,dV\le C\int_{r\varepsilon\leq|z_j|\leq R\varepsilon}|z_j|^{2a_j}\,dV_{z_j}\le C\varepsilon^{2a_j+2}.
\end{aligned}
\]
\end{lemma}

Finally, we extend \cite[Lemma 3.1]{MT08}, after clarifying the distinction between the function spaces.
For a multiplier ideal sheaf $\mathscr I(h)$, we set
\[
A^{p,q}(X,L\otimes\mathscr I(h)):=\Gamma\bigl(X,\mathscr A_X^{p,q}\otimes_{\cal{O}_X}(L\otimes\mathscr I(h))\bigr),
\]
where $\mathscr A_X^{p,q}$ denotes the sheaf of smooth $(p,q)$-forms.
In the horizontal singularity setting (\ref{Setting: horizontal singularity}), the ideal sheaf $\mathscr I(h)$ is locally generated by $s=\prod_{\ell}z_\ell^{\lfloor a_\ell\rfloor}$.
Thus, on such a coordinate neighborhood $U$, we obtain $A^{p,q}(U,L\otimes\mathscr I(h))=s\,A^{p,q}(U,L)$. In particular, these forms are precisely the smooth $L$-valued forms locally divisible by $s$, and $A^{p,q}(V,L\otimes\mathscr I(h))\subset A^{p,q}\cap L^2_{\mathrm{loc}}(V,L,h)$ for any open subset $V\subset X$.
Furthermore, the space $A^{p,q}(X,L\otimes\mathscr I(h))$ is preserved by $\dbar$, and these spaces coincide when restricted to holomorphic sections.
Note that, for any $[\sigma]\in H^0(Y,f_*(\Omega^p_{X/Y}\otimes\lih))$, its admissible representative $\sigma^\sharp$ belongs to $A^{p,0}(X,\lih)$.
In particular, if $[\sigma]\in H^0(Y,\cal{F}^{n-q}_{h,\omega})$, then the proof of Proposition \ref{Proposition: key Prop. MT08, Lemma 4.5} (1) further gives $\delh\sigma^\sharp\in A^{n-q,1}(X,\lih)$.

\begin{lemma}[{cf. \cite[Lemma 3.1]{MT08}}]\label{Lemma: ext MT08, Lemma 3.1}
  Let $\sigma^\sharp\in A^{n-q,0}(X,\lih)$ be an admissible representative of $[\sigma]\in H^0(Y,f_*(\Omega^{n-q}_{X/Y}\otimes\lih))$. Then 
  \begin{itemize}
    \item [(1)] $\dbar\sigma^\sharp=\sum^m_{j=1}\eta^j\wedge dt_j$ for some $\eta^j\in A^{n-q-1,1}(X,\lih)$, 
    \item [(2)] these $\eta^j$ are not unique, but $[\eta^j]\in A^{n-q-1,1}(X/Y,L)$ are well-defined for $\sigma^\sharp$,
    \item [(3)] for any $y$, the forms $\eta^j_y:=\eta^j|_{X_y}$ are $\dbar$-closed and belong to $L^2_{n-q-1,1}(X_y,L_y,h_y)$,
    \item [(4)] there exist $\alpha^j\in A^{n,q+1}(X,\lih)$ such that $\eta^j_y\wedge\omega^{q+1}_y=\dbar \alpha^j_y$ and $\alpha^j_y:=\alpha^j|_{X_y}\in \Dom\dbar\subset L^2_{n,q+1}(X_y,L_y,h_y)$ for any $y\in Y$.
  \end{itemize}
\end{lemma}

\begin{proof}
  (1)--(3) and (4), except for the integrability of $\alpha^j_y$, follow as in \cite[Lemma 3.1]{MT08}. Indeed, since the spaces $A^{p,q}(X,\lih)$ are preserved by $\dbar$, the coefficients of $\eta^j$ take values in $\lih$. 
  Under the Setting (\ref{Setting: horizontal singularity}), locally we write $\varphi=\sum^r_{j=1}a_j\log|z_j|^2+\psi$ and the ideal $\scr{I}(h)$ is generated by $\prod^r_{j=1}z_j^{\gaus{a_j}}$, where $a_j\geq 0$ and $\psi$ is smooth. 
  Since $\eta^j$, $\alpha^j$, and $\dbar\alpha^j$ are smooth forms with values in $\lih$, their restrictions $u_y$ to $X_y$ satisfy locally $|u_y|_{h_y,\omega_y}^2\le C\prod_{j=1}^r|z_j|^{-2(a_j-\gaus{a_j})}$.
  Hence, the result follows. 
\end{proof}


We construct a primitive correction for Lemma \ref{Lemma: ext MT08, Lemma 3.1} that preserves the $\delho$-condition.

\begin{proposition}\label{Proposition: minimal weighted L2 primitive correction}
  Let $\eta^j_0$ and $\alpha^j_0$ be as in Lemma \ref{Lemma: ext MT08, Lemma 3.1}. Let $P_{q+1}$ be the Hilbert orthogonal projection of $\Lxo{q+1}$ onto $\Ker\dbar$ and set $\hajo:=\alpha^j_0-P_{q+1}\alpha^j_0$. 
  Define $\beta^j_0$ by the pointwise Lefschetz isomorphism $\omega^{q+1}_0\wedge\beta^j_0=\hajo$, i.e., $\bjo=(c_{n-q-1}/(q-1)!)\ast_0\hajo$. 
  Then, $\bjo\in L^2_{n-q-1,0}(X_0,L_0,h_0)\cap A^{n-q-1,0}(X_0\setminus D_0,L_0)$ and we obtain $\delho\bjo=0$ and $\omega^{q+1}_0\wedge\dbar\bjo=\omega^{q+1}_0\wedge\eta^j_0$ in the sense of currents on $X_0$. 
\end{proposition}

\begin{proof}
The vector $\hajo$ has minimal norm among the solutions of $\dbar\alpha^j_0=\omega_0^{q+1}\wedge\eta_0^j$ and is orthogonal to $\Ker\dbar_{q+1}$. 
Since $\rom{Im}\,\dbar_q\subset\Ker\dbar_{q+1}$, it follows by the definition of the Hilbert adjoint that $\hajo\in\Dom\dbarado$ and $\dbarado\hajo=0$. 
Clearly, $\bjo\in L^2_{n-q-1,0}(X_0,h_0,h_0)$, and the constructions $\hajo$ and $\bjo$ yield $\omega^{q+1}_0\wedge\dbar\bjo=\dbar\hajo=\omega^{q+1}_0\wedge\eta^j_0$ in the sense of currents. 
Since $\Box_{h_0,\omega_0}\hajo=\big(\dbar\dbarado+\dbarado\dbar\big)\hajo=\dbarado(\omega^{q+1}_0\wedge\eta^j_0)$ is smooth on $X_0\setminus D_0$, interior elliptic regularity implies that $\widehat\alpha_0^j$ and also $\bjo$ are smooth on $X_0\setminus D_0$.

On $X_0\setminus D_0$, the formal adjoint identity $\dbarado=-*_0\partial_{h_0}*_0$ and $\dbarado\hajo=0$ imply $\partial_{h_0}\bjo=0$. 
To extend $\partial_{h_0}\bjo=0$ across $D_0$, we consider local integrability. On a local chart $U\subset X_0$ of $D_0=\{z_1\cdots z_r=0\}$, we write $\varphi_0=\sum^r_{j=1}a_j\log|z_j|^2+\psi_0$, where $a_j\geq0$ and $\psi_0$ is a smooth function.
When $a_j>0$, the Cauchy--Schwarz inequality gives 
\begin{align*}
  \int_U\frac{|\bjo|}{|z_j|}dV\leq\bigg(\int_U|\bjo|^2e^{-\varphi_0}dV\bigg)^{1/2}\bigg(\int_U|z_j|^{-2}e^{\varphi_0}dV\bigg)^{1/2}<+\infty.
\end{align*}
Thus the coefficients of $\partial\varphi_0\wedge \bjo$ are locally integrable, where $\partial\varphi_0=\sum^r_{j=1}a_j\frac{dz_j}{z_j}+\partial\psi_0$. 
Choose a cutoff function \(\chi_\varepsilon\) vanishing near \(D_0\) as in Lemma \ref{Lemma: cutoff near D}. Combining the Cauchy--Schwarz inequality with the estimates in Lemma \ref{Lemma: cutoff near D}, we obtain locally
\[
  \|\partial\chi_\varepsilon\wedge\bjo\|_{L^1}\le C\sum_{\ell=1}^r\varepsilon^{a_\ell}\|\bjo\|_{L^2(h_0;B_{\ell,\varepsilon})}\longrightarrow0.
\]
On the other hand, the identity on $X_0\setminus D_0$ yields $\partial_{h_0}(\chi_\varepsilon\bjo)=\partial\chi_\varepsilon\wedge\bjo$.
Consequently the equality $\partial_{h_0}\bjo=0$ extends to $X_0$ in the sense of currents as $\varepsilon\to0$.
\end{proof}

\section{The Hodge metric and its curvature formula}

When $h$ is smooth, the natural inclusion $f_*(\Omega^p_{X/Y}\otimes L)_y\subset H^0(X_y,\Omega^p_{X_y}\otimes L_y)$ induces on each stalk $f_*(\Omega^p_{X/Y}\otimes L)_y$ 
the canonical $L^2$-pairing given by integration over each fiber $X_y$ with respect to $\omega_y$ and $h_y$.
On the locally free locus, these pairings define a Hermitian metric $g$, called \textit{the Hodge metric}. 
On any locally free subsheaf of $f_*(\Omega^p_{X/Y}\otimes L)$, the induced pairing defines a Hermitian metric wherever the inclusion is fiberwise injective. 
This metric is also called the Hodge metric (see \cite[\S3.3]{MT08}).

In this section, using Corollary~\ref{Corollary: splitting sheaves}, we define a Hodge metric locally on the torsion-free coherent sheaf $R^qf_*(K_{X/Y}\otimes\lih)$. 
In the horizontal singularity setting~(\ref{Setting: horizontal singularity}), we show that the local Hodge metrics define a global smooth Hermitian metric $g$ on $\hdi$, which is locally free by Lemma \ref{Lemma: horizontal base change}. 
Following \cite{MT08}, we then describe the Chern connection of $g$ and derive an optimized curvature formula.

\subsection{Definition of Hodge metrics}

We define a canonical singular Hermitian metric on a torsion-free coherent sheaf $\hdi$, which we call \textit{the Hodge metric}. 

\begin{definition}[local Hodge metric]\label{Definition: Hodge metric}
For every point $y\in Y$, take a local coordinate neighborhood $W\cong\{t=(t_1,\ldots,t_m)\in\bb{C}^m\mid \|t\|<1\}$ centered at $y$, 
and a K\"ahler metric $\omega=\omega_f+cf^*(i\sum_jdt_j\wedge d\overline{t}_j)$ on $f^{-1}(W)$ for a real number $c$.
The choice of $\omega$ gives an injective homomorphism
\begin{align*}
  S_\omega:=\jmath\circ S^q_{f,\omega,h}:\hdi\longrightarrow f_*(\Omega^{n-q}_{X/Y}\otimes\lih)\hookrightarrow f_*(\Omega^{n-q}_{X/Y}\otimes L)
\end{align*}
over $W$ (see Corollary \ref{Corollary: splitting sheaves}). Let $Y^\circ_{q,h}\subset Y$ denote the maximal Zariski open subset where $\hdi$ is locally free.
For every pair of vectors $u_y,v_y\in\hdi_y$ with $y\in W\cap Y^\circ_{q,h}$, we define
\begin{align*}
  g(u_y,v_y)=\frac{c_{n-q}}{q!}\int_{X_y}\omega_f^q\wedge S_\omega(u_y)\wedge h\,\overline{S_\omega(v_y)}\big|_{X_y},
\end{align*}
whenever the fiber integral is finite. In particular, the fiber integral is finite whenever $\scr{I}(h)|_{X_y}=\scr{I}(h_y)$.
Thus $g$ is the pullback of the fiberwise Hodge pairing via $S_\omega$, defined on $W\cap Y^\circ_{q,h}$.
We call $g$ \textit{the Hodge metric on} $\hdi|_W$.
\end{definition}

\begin{remark}\label{Remark: smoothness of Hodge metric g}
As in \cite[Lemma~2.2]{CGP22} and \cite[\S5.1]{WZ25}, under the horizontal singularity setting (\ref{Setting: horizontal singularity}), factoring out the vanishing factors guaranteed by the multiplier ideal sheaf $\scr{I}(h_y)$
from the coefficients of the holomorphic sections leaves locally integrable singular factors of the form $\prod_{j=1}^r |z_j|^{-2\beta_j}, 0\leq\beta_j<1$, which are independent of the base variables $t$, while the remaining factors are smooth.
Thus, derivatives of any order with respect to $t$ admit locally uniform integrable majorants, and differentiation commutes with fiber integration. Consequently, the Hodge metric $g$ is smooth.

Furthermore, as explained in $\S$\ref{subsection: Occurrence horizontal setting}, there exists a Zariski-open subset $Y^\circ\subset Y$ over which the horizontal singularity setting (\ref{Setting: horizontal singularity}) holds.
By Lemma \ref{Lemma: horizontal base change}, we have $Y^\circ\subset Y^\circ_{q,h}$. In particular, $g$ is smooth on $Y^\circ$.
\end{remark}

Under the horizontal singularity setting (\ref{Setting: horizontal singularity}), the homomorphisms $L^q_{f,\omega,h}$ and $S^q_{f,\omega,h}$ depend only on the Kähler class $\kappa_f=[\omega]$, which is a weaker assumption on $\omega$ than that in Lemma \ref{Lemma: independence under pullback}.
This result is analogous to \cite[Lemma 5.2]{MT08}.

\begin{lemma}\label{Lemma: independence in relative Kahler class}
We assume that $h$ satisfies the horizontal singularity setting (\ref{Setting: horizontal singularity}).
Let $\omega_1$ and $\omega_2$ be Kähler forms representing the same relative Kähler class $[\omega_1]=[\omega_2]=\kappa_f\in H^2(X,\mathbb R)$.
Then we have $L^q_{f,\omega_1,h}=L^q_{f,\omega_2,h}$ and $S^q_{f,\omega_1,h}=S^q_{f,\omega_2,h}$.

In particular, the induced Hodge metric $g$ does not depend on the choice of a Kähler form $\omega$ in the relative Kähler class $\{\omega_f|_{X_W}\}\in H^2(X_W,\mathbb R)$ satisfying $\omega|_{X_y}=\omega_f|_{X_y}$ for every $y\in W$.
Consequently, the Hodge metric $g$ defines a global Hermitian metric on $R^qf_*(K_{X/Y}\otimes\lih)$ over $Y$.
\end{lemma}

\begin{proof}
We first prove the independence of $L^q_{f,\omega,h}$.
Since $[\omega_1]=[\omega_2]$, the $\partial\dbar$-lemma on each compact Kähler fiber $X_y$ gives $\omega_{1,y}-\omega_{2,y}=\idd\psi_y$ for some smooth real function $\psi_y$ on $X_y$.
For any relative holomorphic form $[\sigma]\in H^0(X,\Omega^{n-q}_{X/Y}\otimes\lih)$, we have
\begin{align*}
  (\omega_{1,y}^q-\omega_{2,y}^q)\wedge\sigma_y=-i\dbar_z\Big(\partial_z\psi_y\wedge\sum_{k=0}^{q-1}\omega_{1,y}^{q-1-k}\wedge\omega_{2,y}^k\wedge\sigma_y\Big).
\end{align*}
The primitive is locally $L^2$ with respect to $h_y$, since $\sigma_y$ has coefficients in $\scr{I}(h_y)$.
Thus $[\omega_{1,y}^q\wedge\sigma_y]=[\omega_{2,y}^q\wedge\sigma_y]\in H^q(X_y,K_{X_y}\otimes L_y\otimes\scr{I}(h_y))$.
By Lemma \ref{Lemma: horizontal base change}, the natural base change maps are isomorphisms for every $y\in Y$, and $R^qf_*(K_{X/Y}\otimes L\otimes\scr{I}(h))$ is locally free. 
Hence the fiberwise equality implies $L^q_{f,\omega_1,h}=L^q_{f,\omega_2,h}$. 

The independence of $S^q_{f,\omega,h}$ follows from that of $L^q_{f,\omega,h}$ by the same argument as in the proof of Lemma \ref{Lemma: independence under pullback}.
\end{proof}

At this point, over the Zariski-open subset $Y^\circ\subset Y$ introduced in $\S$\ref{subsection: Occurrence horizontal setting}, the holomorphic Hermitian vector bundle $(\hdi,g)|_{Y^\circ}$ admits the Chern connection $D_g=\partial_g+\dbar$ and the Chern curvature $\iO{g}$.
We are interested in the positivity of $\iO{g}$, for which it suffices to consider the following setting.

Let $f:(X,\omega)\longrightarrow Y\subset\bb{C}^m$ be as in local setting (\ref{Setting: local}), and suppose that $(L,h)$ satisfies the horizontal singularity setting (\ref{Setting: horizontal singularity}).
Denote by $\cal{F}:=\cal{F}^{n-q}_{h,\omega}=\rom{Im}\,S^q_{f,\omega,h}$ the image of the homomorphism $S^q_{f,\omega,h}:\hdi \longrightarrow f_*(\Omega^{n-q}_{X/Y}\otimes\lih)$.
By Definition \ref{Definition: Hodge metric}, the canonical pairing $g$ on $f_*(\Omega^{n-q}_{X/Y}\otimes L)$ induces the Hodge metric $g$ on $\hdi$, and this metric is smooth by Remark \ref{Remark: smoothness of Hodge metric g}.
We shall refer to $(\hdi,g)\subset(f_*(\Omega^{n-q}_{X/Y}\otimes L),g)$ as a holomorphic Hermitian subbundle.

\subsection{The Chern connection}

Assume that $h$ satisfies the horizontal singularity setting (\ref{Setting: horizontal singularity}).
As in \cite[Section 5.2]{MT08}, we consider the fiberwise orthogonal projection $P_y:A^{n-q,0}\cap L^{2,loc}_{n-q,0}(X_y,L_y,h_y)\longrightarrow \cal{F}_y$ given by $P_y(u)=\sum_{j=1}^rg_y(u,\sigma_{j,y})\,\sigma_{j,y}$, where $\{\sigma_{j,y}\}_{j=1}^r$ is a $g_y$-orthonormal basis of $\cal{F}_y$.
Since $\cal{F}$ is locally free by Lemma \ref{Lemma: horizontal base change}, the family $\{P_y\}_{y\in Y}$ induces a map
\begin{align*}
  P:A^{n-q,0}\cap L^2_{n-q,0}(X,L,h)\longrightarrow \big\{u\in A^{n-q,0}(Y/Y,L)\big|\, u|_{X_y}\in\cal{F}_y,\forall\,y\in Y\big\}\cong A^0(Y,\cal{F}).
\end{align*}
Then for $[\sigma]\in H^0(Y,\cal{F})$ with $\delh\sigma^\sharp=\sum_j\mu^j\wedge dt_j$ in the sense of currents for some $\mu^j\in A^{n-q,0}\cap L^{2,loc}_{n-q,0}(X,L,h)$ (see Proposition \ref{Proposition: key Prop. MT08, Lemma 4.5}), we define 
\begin{align*}
  \delg[\sigma]:=\sum_j[P(\mu_j)]dt_j \in A^{1,0}(Y,\cal{F}).
\end{align*}

\begin{proposition}[{cf. \cite[Lemmas 5.3--5.4]{MT08}}]\label{Proposition: metric connection}
Assume that $h$ satisfies the horizontal singularity setting (\ref{Setting: horizontal singularity}).
Then the class $[P(\mu^j)]$ is well-defined for $[\sigma]\in H^0(Y,\cal{F})$, and $D_g:=\partial_g+\dbar$ is the Chern connection of the Hodge metric $g$ on $\cal{F}$.
\end{proposition}

\begin{proof}
First, we show that $\partial_g[\sigma]$ is well-defined. It suffices to show that $\mu^j|_{X_y}$ is orthogonal to $\cal{F}_y$ under the condition $[\sigma]=[0]$.
We write $\sigma^\sharp=\sum_j\sigma_j\wedge dt_j$ with some $\sigma_j\in A^{n-q-1,0}\cap L^{2,loc}_{n-q-1,0}(X,L,h)$, then we can take $\mu^j=\delh\sigma_j$. 
For any holomorphic section $s\in\cal{F}_y$, we have
\begin{align*}
  (\mu^j|_{X_y},s)_{h_y}=\frac{c_{n-q}}{q!}\int_{X_y}\omega_y^q\wedge\partial_{h_y}\sigma_j\wedge h_y\overline{s}
    =\frac{c_{n-q}}{q!}\int_{X_y}\partial\left(\omega_y^q\wedge\sigma_j\wedge h_y\overline{s}\right)=0.
\end{align*}
The last equality follows from the cutoff argument in the proof of Proposition \ref{Proposition: key Prop. MT08, Lemma 4.7} (3). 
Indeed, the fiberwise $L^2$-integrability of $\sigma_j|_{X_y}$ under the Setting (\ref{Setting: horizontal singularity}) and of $\mu^j|_{X_y}=\partial_{h_y}\sigma_j$ by Proposition \ref{Proposition: key Prop. MT08, Lemma 4.5} (3), 
together with the vanishing orders of $s$ given by Lemma \ref{Lemma: key vanishing orders}, ensure that the boundary terms tend to zero. Here, Proposition \ref{Proposition: key Prop. MT08, Lemma 4.7} (1) gives $\delhy s=0$. Thus $P_y(\mu^j|_{X_y})=0$ when $[\sigma]=[0]$, proving the well-definedness.

Finally, we show that $D_g$ is compatible with $g$. The Setting (\ref{Setting: horizontal singularity}) ensures that the fiber integrals defining $g$ and the projected coefficients are smooth in the base variables.
Furthermore, differentiation under the fiber integral is justified by the same local integrability and cutoff argument as above. Let $[\sigma],[\tau]\in H^0(Y,\cal{F})$, and we can write $\delh\tau^\sharp=\sum_j\mu^j(\tau)\wedge dt_j$. 
Therefore, the computation in \cite[Lemma 5.4]{MT08} gives
\begin{align*}
  \dbar g([\sigma],[\tau])=(-1)^{n-q}\frac{c_{n-q}}{q!}f_*\Big(\omega^q\!\wedge\!\sigma\!\wedge h\,\overline{\delh\tau^\sharp}\Big)
    \!=\!\sum_jg([\sigma],\overline{[P(\mu^j(\tau))]})\,d\overline{t_j}=g([\sigma],\overline{\partial g[\tau]}),
\end{align*}
where $\sigma|_{X_y}\in\cal{F}_y$. This gives the compatibility with $g$. Together with the Leibniz rule and the fact that the $(0,1)$-part is $\dbar$, the corresponding connection is the Chern connection of $(\cal{F},g)$.
\end{proof}

\subsection{Curvature formula}

Let $Y\subset\mathbb C^m$ be a unit ball centered at $0$, with coordinates $t=(t_1,\ldots,t_m)$. Let $(F,g)$ be a Hermitian holomorphic vector bundle on $Y$ with Chern curvature $\Theta_g=\sum_{j,k}\Theta_{jk}\,dt_j\wedge d\overline{t_k}$, where $\Theta_{jk}\in\rom{End}(Y,F)$.
Then $(F,g)$ is Nakano semi-positive at $t=0$ if and only if $\Theta_g(s)=\sum_{j,k}g_0(\Theta_{jk}s_j,s_k)\ge0$ for any tensor $s=\sum^m_{j=1}\partial/\partial t_j\otimes s_j\in (T_Y\otimes F)_0$. 
On the other hand, for any $\sigma,\tau\in H^0(Y,F)$, we have $\frac{\partial^2}{\partial t_j\partial\overline{t_k}}g(\sigma,\tau)=g((\delg\sigma)^j,(\delg\tau)^k)-g(\Theta_{jk}\sigma,\tau)$, where $\delg\sigma=\sum_j(\delg\sigma)^jdt_j\in A^{1,0}(Y,F)$, etc. 
Hence, if $\sigma$ and $\tau$ are \textit{normal} at $0$ with respect to $g$, i.e., $\delg\sigma=\delg\tau=0$ at $0$, we obtain $\frac{\partial^2}{\partial t_j\partial\overline{t_k}}g(\sigma,\tau)|_{t=0}=-g_0(\Theta_{jk}\sigma|_{t=0},\tau|_{t=0})$. 
Using this fact, we establish Nakano semi-positivity by proving the semi-positivity of $-\idd$ applied to the following $(m-1,m-1)$-form $T_g$.
For $\sigma_1,\ldots,\sigma_m\in H^0(Y,F)$, let $\bm{\sigma}:=(\sigma_1,\ldots,\sigma_m)$ and we define
\[
T_g(\bm{\sigma}):=\sum_{j,k}g(\sigma_j,\sigma_k)\widehat{dt_j\wedge d\overline{t_k}}\in A^{m-1,m-1}(Y,\bb{C}).
\]
In case all $\sigma_j$ are normal at $t=0$, we have $\idd T_g(\bm{\sigma})|_{t=0}=-\sum_{j,k}g_0(\Theta_{jk}\sigma_j|_{t=0},\sigma_k|_{t_0})$.  
Therefore, the Hermitian holomorphic vector bundle $(F,g)$ is Nakano semi-positive at $t=0$, if for any $m$-ple vectors $s_1,\ldots,s_m\in F_0$, there exist extensions $\sigma_j\in H^0(Y,F)$ of $s_j$, all of which are normal at $t=0$ and satisfy $-\idd T_g(\bm{\sigma})|_{t=0}\geq0$ (see \cite[$\S$2]{Ber09}).

\medskip

In what follows, we assume the horizontal singularity setting (\ref{Setting: horizontal singularity}). 
Here, $\cal{F}:=\Imege S^q_{f,\omega,h}$ is a holomorphic vector bundle on $Y$ with a smooth Hodge metric $g$. 
Recall that, for any $[\sigma],[\tau]\in H^0(Y,\mathcal F)$ with representatives $\sigma,\tau\in A^{n-q,0}(X,L\otimes\mathscr I(h))$, we have 
\begin{align*}
  g([\sigma],[\tau])=\frac{c_{n-q}}{q!}f_*\big(\omega_f^q\wedge\sigma\wedge h\overline{\tau}\big).
\end{align*}

\begin{notation}\label{Notation: curvature formula}
  Let $[\sigma_1],\ldots,[\sigma_m]\in H^0(Y,\cal{F})$ with $\sigma_j\in A^{n-q,0}(X,\lih)$. 
  \begin{itemize}
    \item [(1)] Let $N=n-q+m-1$. We set $\wh{\sigma}=\sum_j\sigma_j\wedge\wh{dt_j}\in A^{N,0}(X,\lih)$, then 
    \begin{align*}
      T_g([\bm{\sigma}])=\sum_{j,k}g([\sigma_j],[\sigma_k])\wh{dt_j\wedge d\overline{t_k}}=\frac{c_N}{q!}f_*\big(\omega^q\wedge\wh{\sigma}\wedge h\overline{\wh{\sigma}}\big).
    \end{align*}
    \item [(2)] We write $\delh\sigma_j=\sum_k\mu^k_j\wedge dt_k$, then $\delh\wh{\sigma}=\sum_j\mu^j_j\wedge dt=:\mu\wedge dt$. 
    Here, each $\mu_j^k\in A^{n-q,0}(X,\lih)$ requires $\sigma_j$ to be an admissible representative (see Proposition \ref{Proposition: key Prop. MT08, Lemma 4.5} (1)).
    \item [(3)] We write $\dbar\sigma_j=\sum_k\eta^k_j\wedge dt_k$, then $\dbar\wh{\sigma}=\sum_j\eta^j_j\wedge dt=:\eta\wedge dt$ with $\eta\in A^{n-q-1,1}(X,\lih)$.  
  \end{itemize}
\end{notation}

Let $[\sigma],[\tau]\in H^0(Y,\mathcal F)$ with representatives $\sigma,\tau\in A^{n-q,0}(X,L\otimes\mathscr I(h))$, respectively.
Since $g$ depends only on the restrictions to $X_y$ for each $y\in Y$, we may replace $\sigma$ by $\widetilde{\sigma}=\sigma-\sum_j\varsigma^j\wedge dt_j$, where $\varsigma^j\in L^{2,loc}_{n-q-1,0}(X,L,h)\cap A^{n-q-1,0}(X\setminus D,L)$ need not be smooth on $X$.
Indeed, $\widetilde{\sigma}|_{X_y\setminus D_y}=\sigma|_{X_y\setminus D_y}$. Taking $\widetilde{\tau}$ similarly, we obtain
\[
\begin{aligned}
g([\sigma],[\tau])=\frac{c_{n-q}}{q!}f_*\big(\omega^q\wedge\sigma\wedge h\overline{\tau}\big)=\frac{c_{n-q}}{q!}f_*\big(\omega^q\wedge\widetilde{\sigma} \wedge h\overline{\widetilde{\tau}}\big).
\end{aligned}
\]
Let $\sigma^\circ\in L^{2,loc}_{n-q,0}(X,L,h)\cap A^{n-q,0}(X\setminus D,L)$ satisfy $\sigma^\circ|_{X_y\setminus D_y}=\sigma|_{X_y\setminus D_y}$ for every $y\in Y$. 
Such a form can be constructed as above using correction terms $\varsigma^j$. We call $\sigma^\circ$ \textit{a smooth representative of} $[\sigma]$ \textit{on} $X\setminus D$.
We use the same notation as above for smooth representatives. Writing $\wh{\sigma}^\circ=\sum_j\sigma^\circ_j\wedge\wh{dt_j}\in L^{2,loc}_{N,0}(X,L,h)\cap A^{N,0}(X\setminus D,L)$, we obtain 
\begin{align*}
  T_g([\bm{\sigma}])=\frac{c_N}{q!}f_*\big(\omega^q\wedge\wh{\sigma}\wedge h\overline{\wh{\sigma}}\big)=\frac{c_N}{q!}f_*\big(\omega^q\wedge\wh{\sigma}^\circ\wedge h\overline{\wh{\sigma}^\circ}\big).
\end{align*}

Finally, we derive a curvature formula optimized for the horizontal singularity setting (\ref{Setting: horizontal singularity}). 
To this end, we first introduce a logarithmic cutoff function.

\begin{lemma}[{cf. \cite[Lemma 3.7]{CGP22}}]\label{Lemma: logarithmic cutoff near D}
Let $s_j$ be the canonical section of $\mathcal O_X(D_j)$, and let $|s_j|$ denote its norm with respect to a smooth Hermitian metric.
Choose $\chi\in C^\infty(\mathbb R,[0,1])$ such that $\chi=0$ on $(-\infty,0]$ and $\chi=1$ on $[1,\infty)$. For $0<\varepsilon<1$, set
\[
\chi_\varepsilon=\prod_{j=1}^r\chi_{\varepsilon,j}, \qquad \chi_{\varepsilon,j}=\chi\left(\frac{\log|s_j|-2\log\varepsilon}{|\log\varepsilon|}\right).
\]
Then $\chi_\varepsilon$ vanishes near $D$ and converges pointwise to $1$ on $X\setminus D$.
Put $B_{j,\varepsilon}=\{\varepsilon^2\le |s_j|\le\varepsilon\}$. In local coordinates near $D$, where $D=\{z_1\cdots z_r=0\}$, we have
\begin{align*}
|\partial\chi_{\varepsilon,j}|\le \frac{C}{|z_j||\log\varepsilon|}, \qquad
|\dbar\partial\chi_{\varepsilon,j}|\le C\left(\frac{1}{|z_j|^2|\log\varepsilon|^2}+\frac{1}{|\log\varepsilon|}\right),
\end{align*}
where both derivatives are supported in $B_{j,\varepsilon}\subset\{r\varepsilon^2\le |z_j|\le R\varepsilon\}$ for some constants $r,R>0$ independent of $\varepsilon$.
Furthermore, putting $\Upsilon_j=(|z_j|\log(e/|z_j|))^{-1}$, we obtain
\begin{align*}
|\partial\chi_\varepsilon|\le C\sum_{j=1}^r\mathbf 1_{B_{j,\varepsilon}}\Upsilon_j, \qquad 
|\dbar\partial\chi_\varepsilon|\le C\bigg(\!\Big(\sum_{j=1}^r\mathbf 1_{B_{j,\varepsilon}}\Upsilon_j\Big)^2\!+\sum_{j=1}^r\mathbf 1_{B_{j,\varepsilon}}\Upsilon_j\bigg).
\end{align*}
\end{lemma}

In the horizontal singularity setting~(\ref{Setting: horizontal singularity}), this logarithmic cutoff also yields a curvature formula analogous to that in \cite{MT08}.
Even in this setting, however, the correction terms needed to obtain a representative that is both \textit{primitive} and \textit{strongly normal}, 
as in \cite[Proposition~6.3]{MT08}, need not be smooth across the singular locus $D_0$ (see Remark~\ref{Remark: example of not smooth xi}). 
Guided by \cite[Proposition~4.1]{CGP22}, we therefore establish a version of the curvature formula in \cite[Lemma~5.6]{MT08} and \cite[(4.4)]{Ber09} optimized for applying the correction results of $\S$\ref{Section: Strongly normal and primitive sections}.

\begin{theorem}[{Curvature formula}]\label{Theorem: curvature formula}
  Assume that the horizontal singularity setting (\ref{Setting: horizontal singularity}) holds. For $[\sigma_1],\ldots,[\sigma_m]\in H^0(Y,\cal{F})$, let $\sigma^\circ_j\in L^{2,loc}_{n-q,0}(X,L,h)\cap A^{n-q,0}(X\setminus D,L)$ be a smooth representative of $[\sigma_j]$ on $X\setminus D$. 
  Recall that $T_g([\bm{\sigma}])=f_*\big((c_N/q!)\omega^q\wedge\whs^\circ\wedge h\overline{\whs^\circ}\big)$. We write $\whs^\circ\!=\sum^m_{j=1}\sigma^\circ_j\wedge\wh{dt_j}$, $\delh\whs^\circ\!=\mu^\circ\wedge dt$ and $\dbar\whs^\circ\!=\eta^\circ\wedge dt$ on $X\setminus D$. 
  If $\delh\whs^\circ\!,\, \dbar\whs^\circ\!\in L^2_{loc}(X,L,h)$, then the following curvature formula holds in the sense of currents on $Y$:
  \begin{align*}
    -\idd T_g([\bm{\sigma}])=&\,f_*\Big(\frac{c_N}{q!}\omega^q\wedge(\iO{L,h})_{\rom{ac}}\wedge\whs^\circ\wedge h\overline{\whs^\circ}\Big)\\
    &-\biggl\{\int_{X_t\setminus D_t}\frac{c_{n-q}}{q!}\omega^q_t\wedge\mu^\circ_t\wedge h_t\overline{\mu^\circ_t}+\int_{X_t\setminus D_t}\frac{c_{n-q}}{q!}\omega^q_t\wedge\eta^\circ_t\wedge h_t\overline{\eta^\circ_t}\biggr\}dV_t. 
  \end{align*}

  Furthermore, if each fiber integral on the right-hand side is continuous in $t$, the equality holds for every $t$.
\end{theorem}

\begin{proof}
Choose the logarithmic cutoff $\chi_\varepsilon$ from Lemma \ref{Lemma: logarithmic cutoff near D}, and put $S_\varepsilon=\bigcup_{j=1}^r B_{j,\varepsilon}$, 
which satisfies $\rom{supp}\,\partial\dbar\chi_\varepsilon\!\subset\rom{supp}\,\partial\chi_\varepsilon\!\subset S_\varepsilon$. 
Since $\dbar\whs^\circ\!=\eta^\circ\!\wedge dt$, the corresponding terms vanish for bidegree reasons. Thus, the computation in \cite[Lemma 5.8]{MT08} on $X\setminus D$ gives
\begin{align*}
-\idd& f_*\Big(\chi_\varepsilon \frac{c_N}{q!}\omega^q\wedge\whs^\circ\wedge h\overline{\whs^\circ}\Big)=f_*(\chi_\varepsilon Q)+f_*R_\varepsilon, \quad\text{where} \\
Q&:=\frac{c_N}{q!}\omega^q\wedge(\iO{L,h})_{\rom{ac}}\wedge\whs^\circ\wedge h\overline{\whs^\circ}-\frac{c_{n-q}}{q!}\omega^q\wedge\left(\mu^\circ\wedge h\overline{\mu^\circ}+\eta^\circ\wedge h\overline{\eta^\circ}\right)\wedge dV_t,\\
R_\varepsilon&:=-i\frac{c_N}{q!}\Bigl((-1)^N\partial\chi_\varepsilon\wedge\omega^q\wedge\whs^\circ\wedge h\overline{\partial_h\whs^\circ}-\dbar\chi_\varepsilon\wedge\omega^q\wedge\partial_h\whs^\circ\wedge h\overline{\whs^\circ}\\
&\qquad\qquad\qquad\quad+\dbar\partial\chi_\varepsilon\wedge\omega^q\wedge\whs^\circ\wedge h\overline{\whs^\circ}\Bigr).
\end{align*}

Let $\whs^{\circ(j)}$ denote the component of $\whs^\circ$ containing no $dz_j$, and put $\Upsilon_j\!=\!(|z_j|\log(e/|z_j|))^{-1}$.
Lemma \ref{Lemma: curvature formula} below yields $\Gamma_\sigma\!:=\!|\whs^\circ|_{h,\omega}+\sum_j\Upsilon_j|\whs^{\circ(j)}|_{h,\omega}\!\in\! L^2_{loc}.$
Lemma \ref{Lemma: logarithmic cutoff near D} and the identity $dz_j\wedge\whs^\circ\!=\!dz_j\wedge\whs^{\circ(j)}$ give $|\partial\chi_\varepsilon\wedge\whs^\circ|_{h,\omega}\!\le\! C\mathbf1_{S_\varepsilon}\Gamma_\sigma$.
The same estimate holds for the conjugate terms. Hence, the first two terms in $R_\varepsilon$ are bounded by $C\mathbf1_{S_\varepsilon}\Gamma_\sigma|\partial_h\whs^\circ|_{h,\omega}$.
Furthermore, the product rule gives $\dbar\partial\chi_\varepsilon=\sum_j\dbar\partial\chi_{\varepsilon,j}\prod_{k\ne j}\chi_{\varepsilon,k}+\sum_{j\ne k}\partial\chi_{\varepsilon,j}\wedge\dbar\chi_{\varepsilon,k}\prod_{l\ne j,k}\chi_{\varepsilon,l}$.
Thus, Lemma \ref{Lemma: logarithmic cutoff near D} and the identity $dz_j\wedge\whs^\circ=dz_j\wedge\whs^{\circ(j)}$ show that the terms involving $\partial\chi_{\varepsilon,j}\wedge\dbar\chi_{\varepsilon,k}$ are bounded by
$\Upsilon_j\Upsilon_k|\whs^{\circ(j)}|_{h,\omega}|\whs^{\circ(k)}|_{h,\omega}\!\le\!\big(\sum_j\!\Upsilon_j|\whs^{\circ(j)}|_{h,\omega}\big)^2\!\le\!\Gamma_\sigma^2$.
The terms involving $\dbar\partial\chi_{\varepsilon,j}$ are estimated similarly. Therefore, we have
\begin{align*}
\big|\dbar\partial\chi_\varepsilon\wedge\omega^q\wedge\whs^\circ\wedge h\overline{\whs^\circ}\big|\le C\mathbf1_{S_\varepsilon}\Gamma_\sigma^2\,dV_\omega,
\end{align*}
where $\Gamma_\sigma\in L^2_{loc}$ by Lemma~\ref{Lemma: curvature formula}.
Combining these estimates, we obtain
\begin{align*}
|R_\varepsilon|\le C\mathbf1_{S_\varepsilon}\bigl(\Gamma_\sigma^2+\Gamma_\sigma|\partial_h\whs^\circ|_{h,\omega}\bigr)\,dV_\omega.
\end{align*}

For any $\phi\in C_c^\infty(Y)$, the properness of $f$ ensures that $f^{-1}(\rom{supp}\,\phi)$ is compact, and 
\begin{align*}
|\langle f_*R_\varepsilon,\phi\rangle|=\left|\int_X(f^*\phi)R_\varepsilon\right|\le C_\phi\int_{f^{-1}(\rom{supp}\,\phi)\cap S_\varepsilon}\bigl(\Gamma_\sigma^2+\Gamma_\sigma|\partial_h\whs^\circ|_{h,\omega}\bigr)\,dV_\omega\longrightarrow0.
\end{align*}
Indeed, the integrand belongs to $L^1$ by the Cauchy--Schwarz inequality, since $\Gamma_\sigma,\partial_h\whs^\circ\in L^2_{\mathrm{loc}}$.
Moreover, $\mathbf1_{S_\varepsilon}\to0$ almost everywhere, and hence the convergence follows from the dominated convergence theorem.
Furthermore, the assumptions and the smoothness of the regular part of the local weight imply $(c_N/q!)\,\omega^q\wedge\whs^\circ\!\wedge h\overline{\whs^\circ}, Q\in L^1_{loc}$.
Letting $\varepsilon\to0$, we obtain the desired curvature formula in the senseof currents. 
The pointwise equality follows from the continuity of the fiber integrals in $t$ and the smoothness of $g$.
\end{proof}

\begin{lemma}\label{Lemma: curvature formula}
Let $(U;z_1,\ldots,z_n,t)$ be an adapted coordinate neighborhood of $D$,
where $D\cap U=\{z_1\cdots z_r=0\}$ and $\varphi=\sum_{j=1}^r a_j\log|z_j|^2+\psi$ with $a_j\geq0$ and $\psi$ smooth.
Suppose that $u\in A^{n-q,0}(U\setminus D,L)\cap L^{2,loc}_{n-q,0}(U,L,h)$ satisfies $\delh u, \dbar u\in L^2_{loc}(U,L,h)$.
For each $j\le r$, write $u=dz_j\wedge u'_j+u^{(j)}$, where $u^{(j)}$ contains no $dz_j$-component. Then
\[
\frac{u^{(j)}}{|z_j|\log(e/|z_j|)}\in L^2_{loc}(U,L,h).
\]
If $a_j\notin\mathbb Z$, the logarithmic factor can be omitted.
\end{lemma}

\begin{proof}
Put $s=\prod_{j=1}^r z_j^{\gaus{a_j}}$ and $w=\prod_{j=1}^r|z_j|^{-2(a_j-\gaus{a_j})}$. 
Writing $u=sv$, the assumptions imply that $v, \dbar v$ and $\partial v-\sum_{j=1}^r(a_j-\gaus{a_j})\frac{dz_j}{z_j}\wedge v$ belong to $L^2_{\mathrm{loc}}(w)$, where $L^2(U,w):=\big\{g:U\longrightarrow\bb{C} \text{ is measurable }\big| \int_U|g|^2w\,d\lambda<+\infty\big\}$.
Here, $\dbar v$ extends across $D$ in the distributional sense by a cutoff argument and the Cauchy--Schwarz inequality as in the proof of Proposition \ref{Proposition: minimal weighted L2 primitive correction}.

Since $a_j-\gaus{a_j}<1$, each weight $|z_j|^{-2(a_j-\gaus{a_j})}$ belongs to the Muckenhoupt class $A_2$. As these weights depend on distinct coordinates, $w$ also belongs to $A_2$.
The weighted $L^2$ boundedness of the Beurling--Ahlfors transform \cite[Theorem 0.9]{PV02} gives
\[
\left\|\frac{\partial f}{\partial z_j}\right\|_{L^2(w)}\le C_w\left\|\frac{\partial f}{\dbar z_j}\right\|_{L^2(w)}
\]
for any $f\in C_c^\infty$. Write $v=\sum_Iv_I\,dz_I$. For $U'\Subset U$, choose $\chi\in C_c^\infty(U)$ with $\chi=1$ on $U'$.
Applying the preceding estimate to $\chi v_I$, with the other variables regarded as parameters, and using approximation in weighted $L^2$, we obtain
\[
\begin{aligned}
\left\|\frac{\partial v_I}{\partial z_j}\right\|_{L^2(U',w)}\!\!\!\le\left\|\frac{\partial(\chi v_I)}{\partial z_j}\right\|_{L^2(w)}\!\!\!\le C_w\left\|\frac{\partial(\chi v_I)}{\dbar z_j}\right\|_{L^2(w)}\!\!\!
\le C\biggl(\left\|\frac{\partial v_I}{\dbar z_j}\right\|_{L^2(U,w)}\!\!\!\!\!+\|v_I\|_{L^2(U,w)}\!\biggr)\!<\!+\infty.
\end{aligned}
\]
Consequently, $\nabla_{z_j}v_I=(\partial v_I/\partial x_j,\partial v_I/\partial y_j)\in L^2_{loc}(w)$.

Fix $j\notin I$. If $a_j>\gaus{a_j}$, the connection term $(a_j-\gaus{a_j})\frac{dz_j}{z_j}\wedge v$ shows that $v_I$ has zero trace along $z_j=0$.
The weighted Hardy inequality therefore yields
\[
\int_{U'}\frac{|v_I|^2}{|z_j|^2}\,w\,d\lambda \le C\int_U\bigl(|v_I|^2+|\nabla_{z_j}v_I|^2\bigr)w\,d\lambda<+\infty
\]
for every $U'\Subset U$. If $a_j\in\bb{Z}$, the logarithmic Hardy inequality instead gives
\[
\int_{U'}\frac{|v_I|^2}{|z_j|^2\big(\log(e/|z_j|)\big)^2}w\,d\lambda<+\infty.
\]
These estimates remain valid near intersections of the components of $D$, by applying the one-variable inequalities with the other variables as parameters.
\end{proof}

\section{Strongly normal and primitive sections}\label{Section: Strongly normal and primitive sections}

Throughout this section, we work in the horizontal singularity setting (\ref{Setting: horizontal singularity}). 
We construct sections admitting representatives that are both \textit{primitive} and \textit{strongly normal} at $t=0$. 
The curvature formula of Theorem~\ref{Theorem: curvature formula} applies to these representatives.

\subsection{The strongly normal correction term}

In this subsection, we establish an analogue of \cite[Lemma 6.2]{MT08} in our setting to construct a strongly normal correction term.
For this purpose, we first use Taylor expansions to obtain vanishing orders analogous to those in Lemma \ref{Lemma: key vanishing orders}, even for smooth forms that are not holomorphic.

\begin{lemma}\label{Lemma: logarithmic cutoff for nonholomorphic forms}
Let $\mu_0\in A^{n-q,0}\cap L^2_{n-q,0}(X_0,L_0,h_0)$ satisfy $\delho\mu_0=0$ in the sense of currents on $X_0$.
Then there exists a family of cutoff functions $\{\chi_\varepsilon\}_{\varepsilon>0}$, vanishing near $D_0$ and converging pointwise to $1$ on $X_0\setminus D_0$, such that $\chi_\varepsilon\mu_0\longrightarrow\mu_0$ and $\partial\chi_\varepsilon\wedge\mu_0\longrightarrow0$ in $L^2(X_0,L_0,h_0)$
as $\varepsilon\to0$.
\end{lemma}

\begin{proof}
  By Setting (\ref{Setting: horizontal singularity}), locally write $\varphi_0=\sum_{j=1}^r a_j\log|z_j|^2+\psi(z)$, as in Lemma~\ref{Lemma: key vanishing orders}, and let $\mu_0=\sum_{|I|=n-q}\mu_I\,dz_I$.
The smoothness and local integrability of the coefficients of $\mu_0$ imply $\mu_I=O(|z_\ell|^{\gaus{a_\ell}})$ for every $I,\ell$. 
Thus, the Taylor expansion in $z_\ell,\bar z_\ell$ gives 
\begin{align*}
  \mu_I=\!\sum_{\nu=0}^{\gaus{a_\ell}}c_{\nu,I}(z')z_\ell^\nu\bar z_\ell^{\gaus{a_\ell}-\nu}\!+O(|z_\ell|^{\gaus{a_\ell}+1}), \text{ where } \,
  c_{\nu,I}(z')=\frac{1}{\nu!(\gaus{a_\ell}\!-\nu)!}\frac{\partial^{\gaus{a_\ell}}\mu_I}{\partial z_\ell^\nu\dbar z_\ell^{\gaus{a_\ell}-\nu}}\bigg|_{z_\ell=0}
\end{align*}
and $z'$ denotes the remaining variables. 
Recall that, locally, we have $\partial_{h_0}=\partial_z-\partial\varphi_0\wedge=\partial_z-\big(\sum_{j=1}^ra_j\frac{dz_j}{z_j}+\partial_z\psi\big)\wedge$.
As in the proof of Lemma~\ref{Lemma: key vanishing orders}, comparing the coefficients of $dz_\ell\wedge dz_I$ in $\partial_{h_0}\mu_0=0$, for $\ell\notin I$, yields $(z_\ell\partial_{z_\ell}-a_\ell)\mu_I=O(|z_\ell|^{\gaus{a_\ell}+1})$.
Comparing the terms of degree $\gaus{a_\ell}$, we obtain $(\nu-a_\ell)c_{\nu,I}=0$ for all $0\le\nu\le\gaus{a_\ell}$.
Hence, if $a_\ell>\gaus{a_\ell}$, all $c_{\nu,I}$ vanish, and thus $\mu_I=O(|z_\ell|^{\gaus{a_\ell}+1})$ for $\ell\notin I$.
By smoothness, these estimates extend across the other components of $D_0$.
Let $\delta_{j,I}=1$ if $j\notin I$ and $a_j>\gaus{a_j}$, and $\delta_{j,I}=0$ otherwise.
By smoothness and the preceding vanishing orders, we obtain
\[
|\mu_I|\le C\prod_{j=1}^r|z_j|^{\gaus{a_j}+\delta_{j,I}}.
\]

Choose the logarithmic cutoff $\chi_\varepsilon$ from Lemma \ref{Lemma: logarithmic cutoff near D}. Then, we have $|\partial\chi_{\varepsilon,\ell}|\le\frac{C}{|z_\ell||\log\varepsilon|}$, where $\operatorname{supp}\partial\chi_{\varepsilon,\ell}\subset\{r\varepsilon^2\le |z_\ell|\le R\varepsilon\}$.
Since $dz_\ell\wedge dz_I=0$ for $\ell\in I$, only $\ell\notin I$ contributes to the singular part of $\partial\chi_\varepsilon\wedge\mu_0$, while the smooth part is bounded by $C|\mu_0|/|\log\varepsilon|$.
Thus, Fubini's theorem gives
\[
\begin{aligned}
\|\partial\chi_\varepsilon\wedge\mu_0\|_{L^2(h_0)}^2\le\frac{C}{|\log\varepsilon|^2}\bigg(1+\sum_{I}\sum_{\ell\notin I}\int_{r\varepsilon^2}^{R\varepsilon}\xi^{2(\gaus{a_\ell}+\delta_{\ell,I}-a_\ell)-1}\,d\xi\bigg)
=O(|\log\varepsilon|^{-1})\longrightarrow0.
\end{aligned}
\]
Indeed, each integral is $O(|\log\varepsilon|)$ if $a_\ell=\gaus{a_\ell}$, and $O(\varepsilon^{2(1+\gaus{a_\ell}-a_\ell)})$ otherwise.
\end{proof}

\begin{lemma}\label{Lemma: adjoint vanishing}
Under the assumptions of Lemma \ref{Lemma: logarithmic cutoff for nonholomorphic forms}, set $u=(c_{n-q}/q!)\omega^q_0\wedge\mu_0$ which satisfies $\ast_0u=\mu_0$. 
Then we have $u\in\rom{Dom}\,\dbar^*_{h_0,\omega_0}$ and $\dbar^*_{h_0,\omega_0}u=0$. 
\end{lemma}

\begin{proof}
Let $\chi_\varepsilon$ be the cutoffs constructed in Lemma \ref{Lemma: logarithmic cutoff for nonholomorphic forms}, and put $u_\varepsilon=\chi_\varepsilon u$.
Since $u_\varepsilon$ is smooth and supported away from $D$, the usual formal adjoint formula gives
\begin{align*}
  \dbar^*_{h_0,\omega_0}u_\varepsilon=-\ast_0\delho\ast_0u_\varepsilon=-\ast_0\delho(\chi_\varepsilon\mu_0)=-\ast_0(\partial\chi_\varepsilon\wedge\mu_0).
\end{align*}
Hence, Lemma \ref{Lemma: logarithmic cutoff for nonholomorphic forms} yields $u_\varepsilon\longrightarrow u$ and $\dbar^*_{h_0,\omega_0}u_\varepsilon \longrightarrow0$ in $L^2(X_0,L_0,h_0)$.
Since $\dbar^*_{h_0,\omega_0}$ is closed (see \cite[Chapter VIII, Theorem 1.1]{Dem-book}), the assertion follows.
\end{proof}

Recall that $\scr{I}(h)=\cal{O}_X(-\sum^r_{j=1}\gaus{a_j}D_j)$ is invertible and that $\scr{I}(h)|_{X_y}=\scr{I}(h_y)$, under the Setting (\ref{Setting: horizontal singularity}). 
Since $f$ is smooth, it is flat over $Y$. Moreover, $\Omega^{n-q}_{X/Y}\otimes\lih$ is locally free, and hence flat over $Y$.
To construct strongly normal sections, we further assume that the function $y\longmapsto \dim H^0(X_y,\Omega^{n-q}_{X_y}\otimes L_y\otimes\scr{I}(h_y))$ is constant around $t=0$. 
By Grauert's theorem \cite[Chapter 10, Section 5]{GR84}, this assumption implies that $f_*(\Omega^{n-q}_{X/Y}\otimes\lih)$ is locally free around $t=0$, and its fiber $f_*(\Omega^{n-q}_{X/Y}\otimes\lih)_y$ coincides with $H^0(X_y,\Omega^{n-q}_{X_y}\otimes L_y\otimes\scr{I}(h_y))$.
Recall Corollary \ref{Corollary: splitting sheaves} that $f_*(\Omega^{n-q}_{X/Y}\otimes\lih)=\cal{F}\oplus\cal{K}$, where $\cal{K}:=\cal{K}^{n-q}_{h,\omega}=\Ker L^q_{f,\omega,h}$.

\begin{lemma}[{cf. \cite[Lemma 6.1]{MT08}}]\label{Lemma: normality for pseudoeffective metrics. MT08, Lemma 6.1}
Assume that the function $y\longmapsto \dim H^0(X_y,\Omega^{n-q}_{X_y}$ $\otimes L_y\otimes\!\scr{I}(h_y))$ is constant around $t=0$.
Let $[\sigma]\in H^0(Y,\cal{F})$ with $\partial_h\sigma^\sharp=\sum_j\mu^j\wedge dt_j$, where $\mu_0^j:=\mu^j|_{X_0}\in A^{n-q,0}\cap L^2_{n-q,0}(X_0,L_0,h_0)$ and satisfies $\partial_{h_0}\mu_0^j=0$.
If $\partial_g[\sigma]=0$ at $t=0$, then every $\mu_0^j$ is orthogonal to $H^0(X_0,\Omega_{X_0}^{n-q}\otimes L_0\otimes\scr{I}(h_0))$.
\end{lemma}

\begin{proof}
By Proposition \ref{Proposition: key Prop. MT08, Lemma 4.7}, we have an orthogonal direct sum  $H^0(X_0,\Omega_{X_0}^{n-q}\otimes L_0\otimes\scr{I}(h_0))=\cal{F}_0\oplus \cal{K}_0$ with respect to the inner product $(\cdot,\cdot)_{h_0}$.
Since $\partial_g[\sigma]=0$ at $t=0$, the formula for the Chern connection implies that every $\mu_0^j$ is orthogonal to $\cal{F}_0$.

Fix $j$ and put $u=(c_{n-q}/q!)\,\omega_0^q\wedge\mu_0^j$. By Lemmas \ref{Lemma: logarithmic cutoff for nonholomorphic forms} and \ref{Lemma: adjoint vanishing}, we have $u\in\rom{Dom}\,\dbar^*_{h_0,\omega_0}$ and $\dbar^*_{h_0,\omega_0}u=0$.
The case $q=0$ follows immediately from $\cal{K}_0=0$. We henceforth assume $q>0$. 
Let $\tau_0\in \cal{K}_0$, then Proposition \ref{Proposition: key Prop. MT08, Lemma 4.7} (2) provides $b_0\in A^{n,q-1}\cap L^2_{n,q-1}(X_0,L_0,h_0)$ such that $\dbar b_0=(c_{n-q}/q!)\omega_0^q\wedge\tau_0$. 
%
Consequently,
\begin{align*}
  (\mu_0^j,\tau_0)_{h_0}=\big(u,(c_{n-q}/q!)\,\omega_0^q\wedge\tau_0\big)_{h_0}=(u,\dbar b_0)_{h_0}=(\dbar_{h_0}^*u,b_0)_{h_0}=0.
\end{align*}
Thus $\mu_0^j$ is orthogonal to $\cal{K}_0$. Since it is also orthogonal to $\cal{F}_0$, the direct sum decomposition above gives the assertion.
\end{proof}

\begin{proposition}[{cf. \cite[Lemma 6.2]{MT08}}]\label{Proposition: primitive correction. MT08, Lemma 6.2}
  For $[\sigma]\in H^0(Y,\cal{F})$ with $\delh\sigma^\sharp=\sum_j\mu^j\wedge dt_j$, if $\mu^j_0:=\mu^j|_{X_0}$ is orthogonal to $H^0(X_0,\Omega^{n-q}_{X_0}\otimes L_0\otimes\scr{I}(h_0))$ for some $j$, 
  then there exists $\xi^j_0\in L^2_{n-q-1,0}(X_0,L_0,h_0)$ which is smooth on $X_0\setminus D_0$ and satisfies $\delho\xi^j_0=\mu^j_0$ and $\omega^{q+1}_0\wedge\dbar\xi^j_0=0$, i.e., $\dbar\xi^j_0$ is primitive, on $X_0$ in the sense of currents. 
\end{proposition}

\begin{proof}
  Let $u^j_0:=\omega^q_0\wedge\mu^j_0$. Lemma \ref{Lemma: adjoint vanishing} gives $\ujo\in\Dom\dbar^*_{h_0,\omega_0}$ and $\dbarado\ujo=0$.
  Furthermore, \cite[Theorem 1.3]{QZ26} or $\S$\ref{subsection: Absolute setting} shows that the Hodge $\ast$-operator satisfies $\ast_0\cal{H}^{n,q}(X_0,L_0,h_0)\subset H^0(X_0,\Omega^{n-q}_{X_0}\otimes L_0\otimes\scr{I}(h_0))$. Thus, for any $v\in \cal{H}^{n,q}(X_0,L_0,h_0)$, the assumption and the isometry of the Hodge $\ast$-operator yield
  \begin{align*}
    (\ujo,v)_{h_0,\omega_0}=(\ast_0\ujo,\ast_0v)_{h_0,\omega_0}=q!c_{n-q}^{-1}(\mu^j_0,\ast_0v)_{h_0,\omega_0}=0.
  \end{align*}
  Hence, $\ujo$ is orthogonal to $\cal{H}^{n,q}(X_0,L_0,h_0)$. 
  
  On the other hand, \cite[Proposition 2.18]{Mat22} establishes that $\rom{Im}\,\dbar$ is a closed subspace in $L^2_{n,q}(X_0,L_0,h_0)$. Thus, the orthogonal decomposition \cite{Dem-book} gives
  \begin{align*}
    L^2_{n,q}(X_0,L_0,h_0)=\rom{Im}\,\dbar\oplus \cal{H}^{n,q}(X_0,L_0,h_0)\oplus\rom{Im}\,\dbarado.
  \end{align*}
  Since $\dbarado \ujo=0$ and $\ujo$ is orthogonal to $\calHo$, there is $\cjo\in\! L^2_{n,q+1}(X_0,L_0,h_0)$ such that $\ujo=\dbarado\cjo$. Choosing $\cjo$ to be the solution of minimal $L^2$-norm, we have $\cjo\in(\Ker\dbarado)^{\perp}=\rom{Im}\,\dbar$, and hence $\dbar \cjo=0$. 
  By the pointwise Lefschetz isomorphism, we can write $\cjo=\omega^{q+1}_0\wedge\zeta$ where $\zeta\in L^2_{n-q-1,0}(X_0,L_0,h_0)$. It follows from $\dbar\cjo=0$ that $\omega^{q+1}_0\wedge\dbar\zeta=0$ in the sense of currents.  
  Since $h_0$ is smooth on $X_0\setminus D_0$, the same computation involving the Hodge $\ast$-operator as in \cite[Lemma 6.2]{MT08} yields $\mu^j_0=-i(q+1)\delho\zeta$. 
  Thus, setting $\xi^j_0:=-i(q+1)\zeta$, we obtain the desired two identities on $X_0\setminus D_0$. Note that $\xi^j_0=-i(c_{n-q-1}/q!)\ast_0\cjo$.

  Finally, we extend the identity to all of $X_0$ as in Proposition \ref{Proposition: minimal weighted L2 primitive correction}. Under the Setting (\ref{Setting: horizontal singularity}), we can write locally $\varphi_0=\sum_{\ell=1}^r a_\ell\log|z_\ell|^2+\psi_0$, here $a_\ell\ge0$ and $\psi_0$ is smooth. 
  Since $\xi_0^j\in L^2_{n-q-1,0}(X_0,L_0,h_0)$, the Cauchy--Schwarz inequality implies that $\xi_0^j/z_\ell\in L^1_{\mathrm{loc}}$ whenever $a_\ell>0$. Therefore, $\partial_{h_0}\xi_0^j=\partial\xi_0^j-\partial\varphi_0\wedge\xi_0^j$ is well-defined in the sense of currents.
  Choose a cutoff function \(\chi_\varepsilon\) vanishing near \(D_0\) as in Lemma \ref{Lemma: cutoff near D}. Combining the Cauchy--Schwarz inequality with the estimates in Lemma \ref{Lemma: cutoff near D}, we obtain locally
\[
  \|\partial\chi_\varepsilon\wedge\xi_0^j\|_{L^1}\le C\sum_{\ell=1}^r\varepsilon^{a_\ell}\|\xi_0^j\|_{L^2(B_{\ell,\varepsilon},h_0)}\longrightarrow0.
\]
On the other hand, the identity on $X_0\setminus D_0$ yields $\partial_{h_0}(\chi_\varepsilon\xi_0^j)=\chi_\varepsilon\mu_0^j+\partial\chi_\varepsilon\wedge\xi_0^j$.
Letting $\varepsilon\to0$, we obtain $\partial_{h_0}\xi_0^j=\mu_0^j$ on all of $X_0$ in the sense of currents.
\end{proof}

\begin{remark}\label{Remark: example of not smooth xi}
  Even if $\mu^j_0$ is smooth, a weighted $L^2$-solution $\xi^j_0$ need not be smooth across $D_0$.
  Indeed, consider the local model with $n=1$, $q=0$, $D_0=\{z=0\}$, and $h_0=|z|^{-1}$.
  Let $\chi$ be a smooth cutoff function identically equal to $1$ near the origin and set $\xi^j_0:=\chi|z|$.
  Then $\xi^j_0$ is $L^2$-integrable with respect to $h_0$ and $\partial_{h_0}\xi^j_0=|z|\partial\chi=:\mu^j_0$ in the sense of currents.
  Since $\partial\chi$ is supported away from the origin, $\mu^j_0$ is smooth, whereas $\xi^j_0$ is not even of class $C^1$ at the origin.
  Moreover, $\omega_0\wedge\bar\partial\xi^j_0=0$ holds for degree reasons.
\end{remark}

\subsection{$L^2$-integrability of the $\dbar$-derivatives of the primitive corrections}

In this subsection, we establish the $L^2$-integrability of $\dbar\bjo$ and $\dbar\xi^j_0$ with respect to $h_0$, where $\bjo$ and $\xi^j_0$ are the correction terms constructed in Propositions \ref{Proposition: minimal weighted L2 primitive correction} and \ref{Proposition: primitive correction. MT08, Lemma 6.2}, using the Bochner--Kodaira--Nakano identity.

Under the Setting (\ref{Setting: horizontal singularity}), let $D_0=\sum^r_{j=1}D_j$ on the compact fiber $X_0$. Choose a global smooth Hermitian metric $g_j$ on each $\cal{O}_{X_0}(D_j)$ and its canonical section $s_j$. Then we can write $h_0=h^{sm}_0\prod^r_{j=1}|s_j|^{-2a_j}_{g_j}$, where $h^{sm}_0$ is a smooth Hermite metric on $L_0$. 
For any $\varepsilon>0$, define a smooth Hermite metric $h_{0,\varepsilon}:=h^{sm}_0\prod^r_{j=1}(|s_j|^2_{g_j}+\varepsilon^2)^{-a_j}$, then $h_{0,\varepsilon}\leq h_0$ and $h_{0,\varepsilon}\nearrow h_0$ almost everywhere as $\varepsilon\to0$. Hence, there exists a constant $C_0>0$ independent of $\varepsilon$ such that $\iO{L_0,h_{0,\varepsilon}}\geq-C_0\omega_0$, as follows:
\begin{align*}
  \idd\log(|s_j|^2_{g_j}+\varepsilon^2)=-\frac{|s_j|^2_{g_j}}{|s_j|^2_{g_j}+\varepsilon^2}\iO{D_j,g_j}+\frac{\varepsilon^2|s_j|^2_{g_j}}{(|s_j|^2_{g_j}+\varepsilon^2)^2}i\partial\log|s_j|^2_{g_j}\wedge\dbar\log|s_j|^2_{g_j}.
\end{align*}

\begin{lemma}\label{Lemma: L2-integrability of dbar bjo}
  The form $\bjo$ in Proposition \ref{Proposition: minimal weighted L2 primitive correction} satisfies $\dbar\bjo\in L^2_{n-q-1,1}(X_0,L_0,h_0)$.
\end{lemma}

\begin{proof}
Let $\hajoe$ be the minimal solution in $L^2_{n,q+1}(X_0,L_0,h_{0,\varepsilon},\omega_0)$ of $\dbar\hajoe=\omega_0^{q+1}\wedge\eta_0^j$.
The right side is smooth, and smooth Hodge theory on the compact fiber gives a smooth minimal solution satisfying $\dbar_{h_{0,\varepsilon},\omega_0}^*\hajoe=0$.
Since $h_{0,\varepsilon}\leq h_0$, the form $\hajo$ also solves the same $\dbar$-equation in \(L^2\) with respect to \(h_{0,\varepsilon}\), and hence $\|\hajoe\|_{h_{0,\varepsilon},\omega_0}\leq\|\hajo\|_{h_{0,\varepsilon},\omega_0}\leq\|\hajo\|_{h_0,\omega_0}<+\infty$.
The Bochner--Kodaira--Nakano identity in bidegree $(n,q+1)$ and the curvature bound $\iO{L_0,h_{0,\varepsilon}}\geq-C_0\omega_0$ imply
\begin{align*}
 \|\dbar(*_0\hajoe)\|_{h_{0,\varepsilon},\omega_0}^2
 &\leq \|\dbar\hajoe\|_{h_{0,\varepsilon},\omega_0}^2+\|\dbar_{h_{0,\varepsilon},\omega_0}^*\hajoe\|_{h_{0,\varepsilon},\omega_0}^2+C_1\|\hajoe\|_{h_{0,\varepsilon},\omega_0}^2\\
 &\leq \|\omega_0^{q+1}\wedge\eta_0^j\|_{h_0,\omega_0}^2+C_1\|\hajo\|_{h_0,\omega_0}^2.
\end{align*}
Here $C_1$ is independent of $\varepsilon$ and the operator $\dbar\,\ast_0$ is independent of metrics on $L_0$.
By the uniform $L^2$-bound above and monotone convergence of the metrics, $\hajoe$ has a subsequence converging weakly in every fixed $L^2(h_{0,\delta})$ to a form $\widehat a_{0,\infty}^j\in L^2(h_0)$.
Closedness of the distributional $\dbar$ gives $\dbar\widehat a_{0,\infty}^j=\omega_0^{q+1}\wedge\eta_0^j$. Furthermore, weak lower semicontinuity and letting $\delta\to0$  yield $\|\widehat a_{0,\infty}^j\|_{h_0,\omega_0}\leq\|\hajo\|_{h_0,\omega_0}$.
By uniqueness of the minimal solution, we obtain $\widehat a_{0,\infty}^j=\hajo$. Passing to the weak limit in the above $L^2$-inequality of $\dbar(*_0\hajoe)$, and using the closed distributional operator $\dbar\,\ast_0$, we obtain $\dbar(*_0\hajo)\in L^2(h_0)$.  
Recall that $\bjo=(c_{n-q-1}/(q-1)!)\ast_0\hajo$, which completes the proof.
\end{proof}

Before proving the $L^2$-integrability of $\dbar\xi_0^j$, we recall its construction (see the proof of Proposition \ref{Proposition: primitive correction. MT08, Lemma 6.2}).
Let $[\sigma]\in H^0(Y,\mathcal F)$ with $\partial_h\sigma^\sharp=\sum_j\mu^j\wedge dt_j$, and suppose that $\mu_0^j:=\mu^j|_{X_0}$ is orthogonal to $H^0(X_0,\Omega^{n-q}_{X_0}\otimes\liho)$ for some $j$.
Set $u_0^j:=\omega_0^q\wedge\mu_0^j$. Recall that $\ujo\in\Dom\dbarado\cap\Ker\dbarado$ and is orthogonal to $\calHo$. Let $\gamma_0^j$ be the minimal $L^2$-norm solution of $\ujo=\dbarado\cjo$. 
Then $\gamma_0^j\in\Imege\dbar$, and hence $\dbar\gamma_0^j=0$. The correction term is then given by $\xi_0^j:=-i(c_{n-q-1}/q!)\ast_0\gamma_0^j$.

\begin{lemma}\label{Lemma: L2-integrability of dbar xi}
  The form $\cjo$ constructed above satisfies $\dbar(\ast_0\cjo)\in L^2_{n-q-1,1}(X_0,L_0,h_0)$.
  Hence $\dbar\xi_0^j\in L^2_{n-q-1,1}(X_0,L_0,h_0)$ follows immediately.
\end{lemma}

\begin{proof}
In this lemma, we recover the minimal weighted $L^2$-norm solution $\cjo$ by resolvent regularization and establish the $L^2$-integrability of $\dbar(\ast_0\cjo)$ using the Bochner--Kodaira--Nakano estimate and weak lower semicontinuity.
For $\varepsilon>0$, let $\Tqe$ denote the maximal distributional $\dbar$-operator from $L^2_{n,q}(X_0,L_0,h_{0,\varepsilon})$ to $L^2_{n,q+1}(X_0,L_0,h_{0,\varepsilon})$, and let $\Tqe^*$ be its Hilbert adjoint. For each $\lambda>0$ define
\begin{align*}
 \wel:=(\Tqe^*\Tqe+\lambda\,\rom{id})^{-1}\ujo \quad\text{and}\quad \cel:=\Tqe\wel.
\end{align*}
Since $\Tqe$ is closed and densely defined, the operator $\Tqe^*\Tqe$ is nonnegative and self-adjoint. Hence, the resolvent $(\Tqe^*\Tqe+\lambda\,\rom{id})^{-1}$ exists for any $\lambda>0$ by \cite[Chapter V]{Kat95}. 
It follows that $\Tqe^*\cel=\ujo-\lambda\wel$ and $\dbar\cel=0$, which imply $(\Delta''_{h_{0,\varepsilon},\omega_0}+\lambda)\cel=\dbar \ujo$.
Since $u_0^j$ is smooth, elliptic regularity yields $\cel\in A^{n,q+1}(X_0,L_0)$.

The Bochner--Kodaira--Nakano identity and the inequality $\iO{L_0,h_{0,\varepsilon}}\geq-C_0\omega_0$ show
\begin{align*}
 \|\dbar(*_0\cel)\|_{h_{0,\varepsilon},\omega_0}^2&\leq\|\ujo-\lambda\wel\|_{h_{0,\varepsilon},\omega_0}^2+C_2\|\cel\|_{h_{0,\varepsilon},\omega_0}^2. \tag{$\dagger$}
\end{align*}
Taking the weighted inner product of the resolvent equation with $\wel$ gives
\begin{align*}
  \|\cel\|_{h_{0,\varepsilon},\omega_0}^2+\lambda\|\wel\|_{h_{0,\varepsilon},\omega_0}^2=(\ujo,\wel)_{h_{0,\varepsilon},\omega_0}.
\end{align*}
For fixed $\lambda>0$, this yields uniform bounds in $\varepsilon$ for both $\wel$ and $\cel$, as $h_{0,\varepsilon}\leq h_0$.
In fact, the Cauchy--Schwarz implies $\lambda\|\wel\|_{h_{0,\varepsilon}}\leq\|\ujo\|_{h_{0,\varepsilon}}\leq\|\ujo\|_{h_0}<+\infty$, where $\ujo\in A^{n,q}\cap L^2_{n,q}(X_0,L_0,h_0)$. 
Furthermore, this inequality yields $\lambda\|\cel\|^2_{h_{0,\varepsilon}}\leq\|\ujo\|^2_{h_0}<+\infty$. 
Passing to a common subsequence, we obtain weak limits in $L^2(h_{0,\delta})$ for every $\delta>0$, denoted by $\wl$ and $\cl$. Since $\dbar$ is closed, we have $\cl=\dbar_q\wl$, labeled as $(\ast)$.

For completeness, to verify the adjoint identity in the limit, let $v\in\Dom\dbar_q$.  It also belongs to $\Dom\Tqe$, since $h_{0,\varepsilon}\leq h_0$. The variational equation is
\[
 (\cel,\dbar v)_{h_{0,\varepsilon},\omega_0}+\lambda(\wel,v)_{h_{0,\varepsilon},\omega_0}=(\ujo,v)_{h_{0,\varepsilon},\omega_0}.
\]
By the standard weak convergence argument for increasing weights, we have $(\wel,v)_{h_{0,\varepsilon},\omega_0}$ $\longrightarrow (\wl,v)_{h_0,\omega_0}$. 
The same holds for $(\cel,\dbar v)_{h_{0,\varepsilon},\omega_0}$, yielding $\dbar_q^*\cl=\ujo-\lambda \wl$, labeled as $(\ast\ast)$.
The same convergence argument, together with $(\ast)$ and $(\ast\ast)$, gives 
\begin{align*}
  \|\cel\|_{h_{0,\varepsilon}}^2+\lambda\|\wel\|_{h_{0,\varepsilon}}^2=(\ujo,\wel)_{h_{0,\varepsilon},\omega_0}\longrightarrow (\ujo,\wl)_{h_0,\omega_0}=\|\cl\|^2_{h_0}+\lambda\|\wl\|^2_{h_0}.
\end{align*}
Combining this convergence with weak lower semicontinuity for increasing weights, i.e., $\|\kappa\|^2_{h_0}\!\leq\liminf_{\varepsilon\to0}\|\kappa_\varepsilon\|^2_{h_{0,\varepsilon}}$, we obtain $\|\kappa_\varepsilon\|^2_{h_{0,\varepsilon}}\!\!\longrightarrow\!\|\kappa\|^2_{h_0}$, where $(\kappa,\kappa_\varepsilon)=(\cl,\cel)$ or $(\wl,\wel)$.
Expanding the squared norm yields $\|\ujo-\lambda\wel\|^2_{h_{0,\varepsilon},\omega_0}\!\!\longrightarrow\!\|\ujo-\lambda\wl\|^2_{h_0,\omega_0}$.
Hence, the weak limit of $(\dagger)$ and weak lower semicontinuity yields
\begin{align*}
  \|\dbar(*_0\cl)\|_{h_0,\omega_0}^2 \leq \|\ujo-\lambda \wl\|_{h_0,\omega_0}^2 +C_2\|\cl\|_{h_0,\omega_0}^2. \tag{$\ddagger$}
\end{align*}

The identities $(\ast)$ and $(\ast\ast)$ show that $\cl$ is the Tikhonov-regularized solution of $\dbar_q^*\phi=u_0^j$; see \cite[Section 5.1]{EHN96}.
Indeed, for the Tikhonov functional
\[
F_\lambda(\phi):=\|\dbar_q^*\phi-u_0^j\|_{h_0,\omega_0}^2+\lambda\|\phi\|_{h_0,\omega_0}^2, \quad \phi\in\operatorname{Dom}T_{0,q}^*,
\]
the identities $(\ast)$ and $(\ast\ast)$ gives $F_\lambda(\phi)-F_\lambda(\cl)\!=\!\|\dbar_q^*(\phi-\cl)\|_{h_0,\omega_0}^2\!+\lambda\|\phi-\cl\|_{h_0,\omega_0}^2\!\ge0$.
Thus $\cl$ is the unique minimizer. Comparing it with the minimal weighted $L^2$-norm solution $\cjo$ of $\dbar_q^*\cjo=u_0^j$, we obtain $F_\lambda(\cl)\le F_\lambda(\cjo)=\lambda\|\cjo\|_{h_0,\omega_0}^2$. Consequently,
\[
\begin{aligned}
\|\cl\|_{h_0,\omega_0}\le\|\cjo\|_{h_0,\omega_0}, \quad \|\dbar_q^*\cl-u_0^j\|_{h_0,\omega_0}\le\sqrt{\lambda}\|\cjo\|_{h_0,\omega_0}.
\end{aligned}
\]

Each $\cl$ belongs to $\overline{\operatorname{Im}\dbar_q}$. Every weak limit therefore solves $\dbar_q^*\phi=\ujo$ in this closed subspace, hence equals the unique minimal solution $\cjo$.
The norm bound and weak lower semicontinuity imply strong $L^2(h_0)$ convergence. Letting $\lambda\searrow0$ in $(\ddagger)$ proves 
\[
\|\dbar(\ast_0\cjo)\|^2_{h_0,\omega_0}\leq\|\ujo\|^2_{h_0,\omega_0}+C_2\|\cjo\|^2_{h_0,\omega_0}<+\infty.
\] 
Furthermore, the equality $\xi_0^j:=-i(c_{n-q-1}/q!)\ast_0\gamma_0^j$ implies $\|\dbar\xi^j_0\|^2_{h_0,\omega_0}<+\infty$.
\end{proof}

\subsection{Simultaneously primitive and strongly normal representatives}

Finally, we use the preceding results to extend \cite[Lemma~6.3]{MT08} as follows. 
This extension yields representatives to which the curvature formula of Theorem~\ref{Theorem: curvature formula} applies.

\begin{theorem}\label{Theorem: primitivity and strong normality}
  Under the horizontal singularity setting (\ref{Setting: horizontal singularity}), assume that the function $y\longmapsto\dim H^0(X_y,\Omega^{n-q}_{X_y}\otimes L_y\otimes\scr{I}(h_y))$ is constant around $t=0$. 
  Let $\sigma_0\in\cal{F}_0$ be a vector at $t=0$. Then, after shrinking $Y$ if necessary, there exist a holomorphic section $[\sigma]\in H^0(Y,\cal{F})$ extending $\sigma_0$ and a smooth representative $\sigma^\circ\in L^{2,loc}_{n-q,0}(X,L,h)\cap A^{n-q,0}(X\setminus D,L)$ of $[\sigma]$ on $X\setminus D$, which satisfy the following properties:
  \begin{itemize}
    \item [(1)] $\sigma|_{X_0}=\sigma_0$ and $\sigma^\circ|_{X_y\setminus D_y}=\sigma|_{X_y\setminus D_y}$ for every $y\in Y$, 
    \item [(2)] $\delg[\sigma]=0$ at $t=0$, 
    \item [(3)] $\delh\sigma^\circ\!, \dbar\sigma^\circ\in L^2_{loc}(X,L,h)$.
  \end{itemize}
  We write $\dbar\sigma^\circ=\sum^m_{j=1}\eta^{\circ j}\wedge dt_j$ and $\delh\sigma^\circ=\sum^m_{j=1}\mu^{\circ j}\wedge dt_j$ on $X\setminus D$, where $\eta^{\circ j}$ and $\mu^{\circ j}$ are smooth on $X\setminus D$. 
  Then, for every $j$, the following hold:
  \begin{itemize}
    \item [(4)] $\eta^{\circ j}|_{X_0}\wedge \omega^{q+1}_0=0$ $($primitivity$)$,
    \item [(5)] $\eta^{\circ j}|_{X_0}\in L^2_{n-q-1,1}(X_0,L_0,h_0)$ $(L^2(h_0)$-integrability$)$, 
    \item [(6)] $\mu^{\circ j}|_{X_0}=0$ $($strong normality$)$.  
  \end{itemize}
  Here, the identities $(4)$ and $(6)$ hold in the sense of currents on $X_0$. 
\end{theorem}

\begin{proof}
  Since $\cal F$ is a holomorphic vector bundle with a smooth Hermitian metric $g$, there exists a holomorphic section $[\sigma]\in H^0(Y,\cal{F})$ such that $[\sigma](0)=\sigma|_{X_0}=\sigma_0$ and $\delg[\sigma]=0$ at $t=0$. 
  We take an admissible representative $\sigma^\sharp\in A^{n-q,0}(X,L)$ of $[\sigma]$. By Lemma \ref{Lemma: ext MT08, Lemma 3.1}, we can write $\dbar\sigma^\sharp=\sum^m_{j=1}\eta^j\wedge dt_j$ for some $\eta^j\in A^{n-q-1,1}(X,\lih)$. 
  For each $j$, Proposition \ref{Proposition: minimal weighted L2 primitive correction} gives a primitive correction term $\beta^j_0\in L^2_{n-q-1,0}(X_0,L_0,h_0)\cap A^{n-q-1,0}(X_0\setminus D_0,L_0)$ for $\eta^j_0:=\eta^j|_{X_0}$ such that $\delho\beta^j_0=0$ and $\omega^{q+1}_0\wedge\dbar\beta^j_0=\omega^{q+1}_0\wedge\eta^j_0$ in the sense of currents on $X_0$. 
  By Proposition \ref{Proposition: key Prop. MT08, Lemma 4.5}, we write $\delh\sigma^\sharp=\sum^m_{j=1}\mu^j\wedge dt_j$ for some $\mu^j\in A^{n-q,0}(X,\lih)$, where $\mu^j_0:=\mu^j|_{X_0}\in A^{n-q,0}\cap L^2_{n-q,0}(X_0,L_0,h_0)$ satisfies $\delho\mu^j_0=0$ in the sense of currents.  
  Therefore, the assumption, Lemma \ref{Lemma: normality for pseudoeffective metrics. MT08, Lemma 6.1} and Proposition \ref{Proposition: primitive correction. MT08, Lemma 6.2} give a correction term $\xi^j_0\in L^2_{n-q-1,0}(X_0,L_0,h_0)\cap A^{n-q-1,0}(X_0\setminus D_0,L_0)$ 
  satisfying $\delho\xi^j_0=\mu^j_0$ and $\omega^{q+1}_0\wedge\dbar\xi^j_0=0$ on $X_0$ in the sense of currents. 

  Furthermore, Lemmas \ref{Lemma: L2-integrability of dbar bjo} and \ref{Lemma: L2-integrability of dbar xi} establish $\dbar\beta^j_0, \dbar\xi^j_0\in L^2_{n-q-1,1}(X_0,L_0,h_0)$ for these correction terms. 
  Hence, we first define the restriction of the final correction term to the central fiber by $\varsigma^j_0=\beta^j_0+\xi^j_0$. Then, we have 
  \begin{align*}
    \delho\varsigma^j_0=\mu^j_0 \qquad \text{and} \qquad \omega^{q+1}_0\wedge\dbar\varsigma^j_0=\omega^{q+1}_0\wedge\eta^j_0 \tag{$\spadesuit$}
  \end{align*}
  in the sense of currents, and $\dbar\varsigma^j_0\in L^2_{n-q-1,1}(X_0,L_0,h_0)$. 
  To extend $\varsigma^j_0$, choose finitely many admissible charts $\{(U_\alpha;z_\alpha,t)\}_\alpha$ around $X_0$, extend the local coefficients of $\varsigma^j_0$ constantly in the $t$ variables, and patch them by a smooth partition of unity $\{\rho_\alpha\}_{\alpha}$ subordinate to $\{U_\alpha\}_\alpha$. 
  Denote the resulting forms on a sufficiently small neighborhood of $X_0$ by $\varsigma^j$. These forms restrict to $\varsigma^j_0$ on $X_0$. In every adapted chart the local weights satisfy $e^{-\varphi_\alpha(z_\alpha,t)}\leq Ce^{-\varphi_\alpha(z_\alpha,0)}$ uniformly for small $t$. 
  Thus, $\varsigma^j,\dbar\varsigma^j$ and $\delh\varsigma^j$ are locally in weighted $L^2$ with respect to $h$, and they are smooth on $X\setminus D$. Indeed, the claims for $\varsigma^j$ and $\dbar\varsigma^j$ follow directly. 
  To verify the $L^2(h)$-integrability of $\delh\varsigma^j$, we locally write $\varphi=\sum_j a_j\log|z_j|^2+\psi$, where $a_j\geq0$ and $\psi$ is smooth. Since the extension is constant in $t$ and the singular terms of $\delh$ are independent of $t$, we have 
  \begin{align*}
    \partial_{h_t,z}\,\varsigma^j(z,t)=\partial_{h_0,z}\,\varsigma^j_0(z)+(\partial_z\psi(z,0)-\partial_z\psi(z,t))\wedge\varsigma^j_0(z), \quad \partial_{h,t}\varsigma^j=-\partial_t\psi(z,t)\wedge\varsigma^j.
  \end{align*}
  Since $\delho\varsigma^j_0=\mu^j_0\in L^2_{n-q,0}(X_0,L_0,h_0)$, the $L^2(h_t)$-norm of $\partial_{h_t,z}\,\varsigma^j$ is locally uniformly bounded. 
  By the construction and Fubini's theorem, we obtain $\delh\varsigma^j\in L^{2,loc}_{n-q,0}(X,L,h)$. 

  Finally, we define a smooth representative $\sigma^\circ$ of $[\sigma]$ on $X\setminus D$ by 
  \begin{align*}
    \sigma^\circ:=\sigma^\sharp-\sum^m_{j=1}\varsigma^j\wedge dt_j\in L^{2,loc}_{n-q,0}(X,L,h)\cap A^{n-q,0}(X\setminus D,L).
  \end{align*}
  It follows immediately that $\sigma^\circ|_{X_y\setminus D_y}\!=\sigma^\sharp|_{X_y\setminus D_y}$ for any $y\in Y$ and $\delh\sigma^\circ\!, \dbar\sigma^\circ\!\in\! L^2_{loc}(X,L,h)$. 
  Computing the fiberwise coefficients of $\dbar\sigma^\circ$ and $\delh\sigma^\circ$ and using $(\spadesuit)$, we obtain 
  \begin{align*}
    \eta^{\circ j}|_{X_0}=\eta^j_0-\dbar\varsigma^j_0, \quad \eta^{\circ j}|_{X_0}\wedge\omega^{q+1}_0=0, \quad \text{and}\quad \mu^{\circ j}|_{X_0}=\mu^j_0-\delho\varsigma^j_0=0.
  \end{align*} 
  The $L^2(h_0)$-integrability of $\eta^{\circ j}|_{X_0}$ follows from $\dbar\varsigma^j_0\in L^2_{n-q-1,1}(X_0,L_0,h_0)$.
\end{proof}

\section{Nakano semi-positivity}

\subsection{Smooth Nakano semi-positivity under the horizontal singularity setting}

Let $f:(X,\omega)\longrightarrow Y\subset\bb{C}^m$ be as in the Setting (\ref{Setting: local}), and suppose that $(L,h)$ satisfies the Setting (\ref{Setting: horizontal singularity}). 
Recall that the Hodge $\ast$-operator induces an injective homomorphism $S^q_{f,\omega,h}:R^qf_*(K_{X/Y}\otimes\lih)\longrightarrow f_*(\Omega^{n-q}_{X/Y}\otimes\lih)$.
By Lemma \ref{Lemma: horizontal base change}, its image $\cal{F}:=\Imege S^q_{f,\omega,h}$ is locally free and carries a smooth Hodge metric $g$ (see Remark \ref{Remark: smoothness of Hodge metric g}).    
%

\begin{theorem}\label{Theorem: smooth Nakano semi-positivity}
  Under the horizontal singularity setting (\ref{Setting: horizontal singularity}), the smooth Hodge metric $g$ on $\cal{F}$ is Nakano semi-positive on $Y$, and hence the induced smooth Hodge metric $g$ on $R^qf_*(K_{X/Y}\otimes\lih)$ is also Nakano semi-positive on $Y$.
\end{theorem}

\begin{proof}
Since $g$ is a smooth Hermitian metric of $\cal{F}$, to show the Nakano semi-positivity, it is enough to show it on the complement of an analytic subset of $Y$. 
By Grauert's theorem \cite[Chapter 10, Section 5]{GR84}, there exists an analytic subset $Z\subset Y$ such that the function $y\longmapsto\dim H^0(X_y,\Omega^{n-q}_{X_y}\otimes L_y\otimes\scr{I}(h_y))$ is constant on $Y\setminus Z$. 
Let $y\in Y\subset Z$ be an arbitrary point, and choose local coordinates $t=(t_1,\ldots,t_m)$ centered at $y$. 
Let $\sigma_{1,0},\ldots,\sigma_{m,0}\in\cal F_0$ be arbitrary.
Applying Theorem \ref{Theorem: primitivity and strong normality} to each $\sigma_{j,0}$, we obtain holomorphic sections $[\sigma_j]$ and their smooth representatives $\sigma_j^\circ$ of $[\sigma_j]$ on $X\setminus D$.
Recall that these representatives $\sigma_j^\circ$ are constructed by extending the correction terms $\varsigma_{j,0}^k$ constantly in the base directions in adapted coordinates and gluing them by a smooth partition of unity. 
In particular, $\sigma_j^\circ:=\sigma_j^\sharp-\sum_{k=1}^m\varsigma_j^k\wedge dt_k$.

We first establish the continuity of the fiber integrals appearing in the curvature formula. Since the local weight is given by $\varphi(z,t)=\sum_j a_j\log|z_j|^2+\psi(z,t)$, we have
\begin{align*}
e^{-\varphi(z,t)}=e^{-\varphi(z,0)}e^{\psi(z,0)-\psi(z,t)}\quad\text{and}\quad C^{-1}e^{-\varphi(z,0)}\le e^{-\varphi(z,t)}\le Ce^{-\varphi(z,0)}
\end{align*}
uniformly in $t$. Furthermore, the singular part of $\partial_h$, i.e., $\sum_ja_jdz_j/z_j$, is independent of $t$.
Consequently, the local terms arising from $\varsigma^k_j,\delh\varsigma^k_j$ and $\dbar\varsigma^k_j$ 
are finite sums of fixed weighted $L^2$ coefficients multiplied by functions smooth in $t$.

Let $\gamma$ denote a general form occurring in these fiber integrals of curvature formula in Theorem \ref{Theorem: curvature formula}, and set $\gamma_t=\gamma|_{X_t}$.
Approximate the fixed $L^2$ coefficients $\gamma_0$ on $X_0$ by smooth coefficients supported away from $D_0$, and extend them constantly in the base directions to obtain smooth forms $\gamma_{\nu}$ supported away from $D$.
The uniform equivalence of the weights and the smooth partition of unity give, on a sufficiently small $Y'\Subset Y$,
\begin{align*}
\delta_\nu:=\sup_{t\in Y'}\|\gamma_t-\gamma_{\nu,t}\|_{L^2(h_t,\omega_t)}\longrightarrow0.
\end{align*}
Put $I(t)=\|\gamma_t\|_{L^2(h_t,\omega_t)}^2$ and $I_\nu(t)=\|\gamma_{\nu,t}\|_{L^2(h_t,\omega_t)}^2$. Each $I_\nu$ is continuous, since it is a fiber integral of smooth approximating forms. Then, from the inequality 
\begin{align*}
|I(t)-I_\nu(t)|&\le\bigl(\|\gamma_t\|_{L^2(h_t,\omega_t)}+\|\gamma_{\nu,t}\|_{L^2(h_t,\omega_t)}\bigr)\|\gamma_t-\gamma_{\nu,t}\|_{L^2(h_t,\omega_t)}\\
&\le\Big(2\sup_{t\in Y'}\|\gamma_t\|_{L^2(h_t,\omega_t)}+\delta_\nu\Big)\delta_\nu,
\end{align*}
where $\sup_{t\in Y'}\|\gamma_t\|_{L^2(h_t,\omega_t)}<+\infty$, we obtain $\sup_{t\in Y'}|I(t)-I_\nu(t)|\longrightarrow0$.
Hence $I$ is continuous as the uniform limit of the continuous functions $I_\nu$.

The same argument, together with the Cauchy--Schwarz inequality, applies to the fiber integrals in the curvature formula, which are finite sums of weighted $L^2$ pairings with smooth uniformly bounded coefficients.
The terms arising from the original admissible representatives are treated similarly using their local holomorphic expressions and the uniform weighted $L^2$ estimates under the horizontal singularity setting (\ref{Setting: horizontal singularity}).
Consequently, the curvature formula of Theorem~\ref{Theorem: curvature formula} can be evaluated at $t=0$.

Furthermore, Theorem \ref{Theorem: primitivity and strong normality} gives $\delg[\sigma_j]=0$ at $t=0$ and $\mu_j^{\circ k}|_{X_0}=0, \eta_j^{\circ k}|_{X_0}\wedge\omega_0^{q+1}=0$ in the sense of currents. 
It follows that $\mu_0^\circ:=\mu^\circ|_{X_0}=0$ and $\eta_0^\circ:=\eta^\circ|_{X_0}$ is primitive. The Hodge identity in \cite[Lemma 2.1]{MT08} then yields
\begin{align*}
-\int_{X_0\setminus D_0}\frac{c_{n-q}}{q!}\omega_0^q\wedge\eta_0^\circ\wedge h_0\overline{\eta_0^\circ}=\|\eta_0^\circ\|_{h_0,\omega_0}^2\ge0.
\end{align*}
Finally, applying Theorem \ref{Theorem: curvature formula} at $t=0$, we obtain
\begin{align*}
-\idd T_g([\bm{\sigma}])\big|_{t=0}=f_*\Big(\frac{c_N}{q!}\omega^q\wedge(\iO{L,h})_{\rom{ac}}\wedge\whs^\circ\wedge h\overline{\whs^\circ}\Big)\Big|_{t=0}+\|\eta_0^\circ\|_{h_0,\omega_0}^2\,dV_0\ge0.
\end{align*}
As $\sigma_{1,0},\ldots,\sigma_{m,0}$ are arbitrary, the Hodge metric $g$ is Nakano semi-positive at $t=0$. 
\end{proof}

\subsection{Definition of singular Hermitian metrics and its positivity}

We consider singular Hermitian metrics on holomorphic vector bundles and torsion-free sheaves in general.
Let $X$ be a complex manifold of dimension $n$. 

\begin{definition}[{\cite{HPS18,PT18,DNWZ23,IMW24}}]
  Let $E$ be a holomorphic vector bundle on $X$ and $h$ be a singular Hermite metric on $E$ that is finite and positive definite almost everywhere (see \cite{HPS18} for the precise definition). 
  The metric $h$ is said to be
  \begin{itemize}
    \item \textit{Griffiths semi}-\textit{negative}, if for any open subset $U\subset X$ and any $s \in H^0(U,E)$, the function $\log|s|^2_h$ on $U$ is plurisubharmonic. 
    \item \textit{Griffiths semi}-\textit{positive}, if the dual metric $h^*$ is Griffiths semi-negative. 
    \item \textit{Nakano semi}-\textit{positive}, if the following holds: for data consisting of 
    \begin{itemize}
      \item a Stein coordinate neighbourhood $S\subset X$; 
      \item a \kah metric $\omega_S$;
      \item a smooth function $\psi$ on $S$ such that $\idd\psi>0$; 
      \item a positive integer $q$ with $1\le q\le n$; 
      \item a $\dbar$-closed $g\in L^2_{n,q}(S,E;\omega_S,he^{-\psi})$ satisfying 
      \[
        \int_S\lara{B^{-1}_{\psi,\omega_S}g}{g}_{h,\omega_S}e^{-\psi}dV_{\omega_S}<+\infty,
      \]
      there exists $u\in L^2_{n,q-1}(S,E;\omega_S,he^{-\psi})$ such that 
      \begin{align*}
        \qquad\qquad\dbar u=g \quad\text{and}\quad \int_S|u|^2_{h,\omega_S}e^{-\psi}dV_{\omega_S}\leq\int_S\lara{B^{-1}_{\psi,\omega_S}g}{g}_{h,\omega_S}e^{-\psi}dV_{\omega_S},
      \end{align*}
      where $B_{\psi,\omega_S}=[\idd\psi\otimes\rom{id}_E,\Lambda_{\omega_S}]$.
    \end{itemize}
  \end{itemize}
\end{definition}

\begin{definition}[{\cite[Definition 3.3]{Wat26}}]\label{Definition: uniformly positive definite}
    A singular Hermitian metric $h$ on $E$ is said to be \textit{locally uniformly positive definite on} $X$ if, for every point of $X$, there exist a neighborhood $U\subset X$ of the point, a positive constant $C_U>0$, 
    and a smooth Hermitian metric $h_0$ on $E|_U$ such that $|\bullet|^2_h\geq C_U|\bullet|^2_{h_0}>0$ almost everywhere on $U$.
    In other words, the norm $|\bullet|^2_{h^*}$ is uniformly bounded from above almost everywhere.
\end{definition}

Clearly, Griffiths semi-positivity implies local uniform positive definiteness.

\begin{definition}[{\cite{HPS18,IMW24}}]
  Let $\cal{E}$ be a torsion-free sheaf on $X$. Denote by $X_{\cal{E}}\subset X$ the maximal analytic Zariski open subset on which $\cal{E}$ is locally free. Note that $\rom{codim}\,(X\setminus X_{\cal{E}})\geq2$. 
  A \textit{singular Hermitian metric on a torsion}-\textit{free sheaf} $\cal{E}$ is defined as a singular Hermitian metric $h$ on the vector bundle $E:=\cal{E}|_{X_{\cal{E}}}$.  
  \begin{itemize}
    \item a singular Hermitian metric $h$ on $\cal{E}$ is said to be \textit{Griffiths} (resp. \textit{Nakano}) \textit{semi}-\textit{positive} if the metric $h$ on $E$ is Griffiths (resp. Nakano) semi-positive. 
    \item a singular Hermitian metric $h$ on $\cal{E}$ has \textit{minimal extension property}, if there is a nowhere dense closed analytic subset $Z\subset X$ with the following two properties: 
    \begin{itemize}
      \item $\cal{E}$ is locally free on $X\setminus Z$, or equivalently, $X\setminus Z\subseteq X_{\cal{E}}$. 
      \item for every embedding $\iota:B\hookrightarrow X$ of the open unit ball $B\subset\bb{C}^n$ with $x=\iota(0)\in X\setminus Z$, and every vector $v\in E_x=\cal{E}_x$ with $|v|_h(x)=1$, 
      there is a holomorphic section $s\in H^0(B,\iota^*\cal{F})$ such that 
      \begin{align*}
        s(0)=v \quad\text{and}\quad\frac{1}{\rom{Vol}(B)}\int_B|s|^2_hdV_B\leq1.
      \end{align*}
    \end{itemize}
  \end{itemize}
\end{definition}

The minimal extension property originates in the sharp forms of the Ohsawa--Takegoshi extension theorem \cite{OT87} due to B{\l}ocki~\cite{Blo13} and Guan and Zhou~\cite{GZ15}. 
It can be viewed as intermediate between Griffiths semi-positivity and Nakano semi-positivity.

\subsection{Griffiths/Nakano semi-positivity of the singular Hodge metric}

In this subsection, we prove \ref{Main Theorem}. By hypothesis, $h$ has analytic singularities modulo smooth functions.
%
As in $\S$\ref{subsection: Occurrence horizontal setting}, take a log resolution $\pi:\tl{X}\longrightarrow X$ adapted to the singularities of $h$.
Then $\pi^*h$ has divisorial singularities modulo smooth functions along a simple normal crossing divisor $\tl{D}$.
Furthermore, there exists a suitable analytic Zariski open subset $Y^\circ\subset Y$ such that the induced morphism $\tl{f}:=f\circ\pi:\tl{f}^{-1}(Y^\circ)\longrightarrow Y^\circ$ is smooth and $\tl{D}$ has relative simple normal crossings over $Y^\circ$.
Consequently, the horizontal singularity setting (\ref{Setting: horizontal singularity}) is satisfied over $Y^\circ$.
%
%
Recall 
that Grauert--Riemenschneider-type vanishing and the Leray spectral sequence yield the isomorphism 
\begin{align*}
  R^qf_*(K_{X/Y}\otimes\lih)\cong R^q\tl{f}_*(K_{\tl{X}/Y}\otimes\pi^*L\otimes\scr{I}(\pi^*h))
\end{align*}
for every $q\geq0$. For simplicity, we write $(\pi^*L,\pi^*h)=(\tl{L},\tl{h})$ and consider the Hodge metric $g$ on $R^qf_*(K_{X/Y}\otimes\lih)$, which is induced by $\tl{f}$.
However, since the construction of the Hodge metric requires the smoothness of the morphism, $g$ is initially defined only on $Y^\circ$, but is globally well-defined and smooth on $Y^\circ$ (see Lemma~\ref{Lemma: independence in relative Kahler class} and Remark \ref{Remark: smoothness of Hodge metric g}). 
Furthermore, by Theorem~\ref{Theorem: smooth Nakano semi-positivity}, the smooth Hodge metric $g$ on $Y^\circ$ is Nakano semi-positive.
Let $Y^\circ_{q,h}$ denote the maximal Zariski open subset of $Y$ on which $R^qf_*(K_{X/Y}\otimes\lih)$ is locally free. We obtain $Y^\circ\subset Y^\circ_{q,h}$ by Lemma~\ref{Lemma: horizontal base change}. 

Set $\Delta:=Y\setminus Y^\circ$, which is an analytic subset of $Y$ and satisfies $Y^\circ_{q,h}\setminus\Delta=Y^\circ$.
We investigate whether the Hodge metric $g$ extends across $\Delta$ to $Y^\circ_{q,h}$ while preserving its positivity in an appropriate sense. 
Our setting is similar to that in \cite{Tak24}. Note that, after taking a further log resolution if necessary, $\tl{f}^{-1}(\Delta)$ has simple normal crossings over generic points of $\Delta$.
To establish Griffiths semi-positivity, however, we take a different approach from that of \cite{Tak24}. 
The argument there relies on results of \cite{MT09}, and it is not immediately clear how to adapt these results to the pseudo-effective case.

We first use the following lemma to show that the Hodge metric $g$ is locally uniformly positive definite on $Y^\circ$, and then use this fact to establish the positivity properties of arbitrary extensions of $g$.

\begin{lemma}\label{Lemma: L2 extension of Hodge sections}
Let $B\subset Y$ be a small ball and let $u\in H^0(B\setminus\Delta,R^qf_*(K_{X/Y}\otimes\lih))$. 
If $\int_{B\setminus\Delta}|u|_g^2\,dV<+\infty$, then $u$ extends uniquely to $\widetilde u\in H^0(B,\hdi)$.
\end{lemma}

\begin{proof}
Using the results of \cite{QZ26}, we adapt the proof of \cite[Lemma~3.4]{Tak24} to the pseudo-effective setting.
By the isomorphism 
\[
R^qf_*(K_{X/Y}\otimes\lih)\cong R^q\tl{f}_*(K_{\tl{X}/Y}\otimes\tl{L}\otimes\scr{I}(\tl{h})),
\] 
we work with the K\"ahler morphism $\widetilde f:\widetilde X\longrightarrow Y$.
After shrinking $B$ if necessary, we can choose a K\"ahler form $\omega$ on $\widetilde X_B:=\widetilde f^{-1}(B)$. Put $B^\circ:=B\setminus\Delta$, and let $t=(t_1,\ldots,t_m)$ be coordinates on $B$.
Here, for the relative K\"ahler form $\widetilde\omega_f$ used to construct $g$, we take $\omega=\widetilde\omega_f+c\widetilde f^*(i\sum_{j=1}^mdt_j\wedge d\overline t_j)$.
By Lemma~\ref{Lemma: independence under pullback}, this choice does not change $g$.

First, choose an analytic hypersurface $\Gamma\subset B$ containing $\Delta\cap B$, and put $V:=B\setminus\Gamma$.
Since $V$ is Stein and $V\subset B^\circ$, the restriction $\widetilde f:\widetilde X_V\longrightarrow V$ is a proper map to a Stein space, where $\widetilde X_V:=\widetilde f^{-1}(V)$.
Consequently, $\widetilde X_V$ is holomorphically convex, and hence \cite{QZ26} applies. In particular, Theorem~\ref{Theorem: QZ26} and Corollary~\ref{Corollary: relative splitting} are applicable.

By Corollary~\ref{Corollary: splitting sheaves}, we obtain $\sigma:=S^q_{\widetilde f,\omega,\widetilde h}(u|_V)\in H^0(\widetilde X_V,\Omega_{\widetilde X/V}^{n-q}\otimes\widetilde L\otimes\mathscr I(\widetilde h))$.
Define the holomorphic form $G_V:=\sigma\wedge\widetilde f^*dt\in H^0(\widetilde X_V,\Omega_{\widetilde X}^{n+m-q}\otimes\widetilde L\otimes\mathscr I(\widetilde h))$.
By the definition of the Hodge metric, Fubini's theorem, and the inner product identity for the Hodge star operator \cite[Corollary 2.2]{MT08}, there exists a constant $C>0$ such that
\begin{align*}
  \int_{\widetilde X_V}|G_V|_{\widetilde h,\omega}^2\,dV_\omega\le C\int_V|u|_g^2\,dV\le C\int_{B^\circ}|u|_g^2\,dV<+\infty.
\end{align*}

Since the local weights of $\widetilde h$ are plurisubharmonic and locally bounded above, $G_V$ also belongs to $L^2_{loc}$ in the ordinary sense.
By the $\dbar$-extension property (see \cite[Lemma 6.8]{Dem82}), $G_V$ extends across the analytic subset $\widetilde f^{-1}(\Gamma)$ to a holomorphic form $G$ on $\widetilde X_B$.
Thus, we obtain $G\in H^0(\widetilde X_B,\Omega_{\widetilde X}^{n+m-q}\otimes\widetilde L\otimes\mathscr I(\widetilde h))$.
It is important to note that the decomposition into relative differential forms given by Theorem~\ref{Theorem: QZ26}~(3) is used only over $V$, where $\widetilde f$ is smooth.
In what follows, we directly use the holomorphic form $G$ extended to $\widetilde X_B$, without requiring a relative decomposition over the singular fibers.

Therefore, Theorem~\ref{Theorem: QZ26}~(2) or \cite[Theorem 1.1]{QZ26} gives
\begin{align*}
  [U]:=\frac{c_{N-q}}{q!}L^q_{\omega,\widetilde h}(G)=\Big[\frac{c_{N-q}}{q!}\omega^q\wedge G\Big]\in H^q(\widetilde X_B,K_{\widetilde X}\otimes\widetilde L\otimes\mathscr I(\widetilde h)),
\end{align*}
where $\widetilde X_B$ is holomorphically convex. This construction does not require the smoothness of $\widetilde f$. 
The Steinness of $B$ implies that $H^q(\widetilde X_B,K_{\widetilde X}\otimes\widetilde L\otimes\mathscr I(\widetilde h))\simeq H^0(B,\hdi)$.
Thus, $[U]$ determines a section $\widetilde u\in H^0(B,\hdi)$.

Since $\widetilde f$ is smooth over $V$, Corollary~\ref{Corollary: relative splitting} applies, and we obtain
\begin{align*}
  [U]|_{\widetilde X_V}=\Big[\frac{c_{n+m-q}}{q!}\omega^q\wedge S^q_{\widetilde f,\omega,\widetilde h}(u|_V)\wedge\widetilde f^*dt\Big]=u|_V\,dt.
\end{align*}
Hence $\widetilde u|_V=u|_V$. Since $V$ is dense in $B^\circ$ and both $\widetilde u|_{B^\circ}$ and $u$ are holomorphic on $B^\circ$, the identity theorem implies $\widetilde u|_{B^\circ}=u$.
Finally, the sheaf $\hdi$ is torsion-free (see \cite[Corollary 1.5]{Mat22}), and hence the extension is unique.
\end{proof}

\begin{proposition}\label{Proposition: local uniform positive definiteness}
Extend $g$ to a singular Hermitian metric $\breve{g}$ on $\hdi$ by assigning arbitrary measurable Hermitian forms on $Y^\circ_{q,h}\setminus Y^\circ$.
Then $\breve{g}$ is locally uniformly positive definite on $Y^\circ_{q,h}$.
More precisely, for every $y_0\in Y^\circ_{q,h}$, there exist a neighborhood $W\subset Y^\circ_{q,h}$ of $y_0$ and a constant $C>0$ such that $\breve{g}\ge Cg_0$ almost everywhere on $W$,
where $g_0$ is any smooth Hermitian metric on $\hdi|_{Y^\circ_{q,h}}$.

Furthermore, the singular Hermitian metric $\breve{g}$ is Nakano semi-positive and has the minimal extension property.
\end{proposition}

\begin{proof}
Fix any $y_0\in Y^\circ_{q,h}$ and choose a sufficiently small coordinate ball $B\Subset Y^\circ_{q,h}$ centered at $y_0$ such that $\hdi|_B$ is trivial. Put $B^\circ:=B\setminus\Delta$.

Following the argument of \cite[Proposition 23.3]{HPS18}, consider the Hilbert space
\[
\mathcal A_g^2(B^\circ):=\left\{u\in H^0(B^\circ,\hdi)\ \middle|\ \int_{B^\circ}|u|_g^2\,dV<+\infty \right\}.
\]
By Lemma~\ref{Lemma: L2 extension of Hodge sections}, every $u\in\mathcal A_g^2(B^\circ)$ extends uniquely to a holomorphic section $\widetilde u$ on $B$. Hence we obtain a linear map
\[
T:\mathcal A_g^2(B^\circ)\longrightarrow H^0(B,\hdi): u\longmapsto\widetilde u.
\]
The graph of $T$ is closed, and hence $T$ is continuous by the closed graph theorem. Indeed, suppose that $u_\nu\longrightarrow u$ in $L^2(B^\circ,g)$ and $T(u_\nu)=\widetilde u_\nu\longrightarrow w$ locally uniformly on $B$.
Since $g$ is smooth on $B^\circ$, convergence in $L^2(g)$ implies local uniform convergence on $B^\circ$. Consequently, $w|_{B^\circ}=u$.
The uniqueness of the holomorphic extension then gives $T(u)=w$.
For any $B'\Subset B$, the seminorm $p_{B'}(w):=\sup_{y\in\overline{B'}}|w(y)|_{g_0}$ is continuous on $H^0(B,R^qf_* (K_{X/Y}\otimes\lih))$.
Therefore, the continuity and linearity of $T$ imply that there exists a constant $C_1>0$ such that
\[
\sup_{y\in B'}|\widetilde u(y)|_{g_0}^2 \le C_1\int_{B^\circ}|u|_g^2\,dV.
\]

On the other hand, $g$ is smooth Nakano semi-positive on $B^\circ$. 
Thus, the optimal $L^2$ extension theorem \cite[Theorem 2.13]{SY23} (cf. \cite[Proposition 3.6]{Tak24}) implies that, for every $y\in B'\cap B^\circ$ and $v\in\hdi_y$, there exists $u\in\mathcal A_g^2(B^\circ)$ such that
\[
u(y)=v\qquad\text{and}\qquad \int_{B^\circ}|u|_g^2\,dV \le C_2|v|_{g_y}^2,
\]
where $C_2>0$ is independent of $y$ and $v$. Combining these two estimates, we obtain $|v|_{g_0}^2 =|\widetilde u(y)|_{g_0}^2\le C_1C_2|v|_{g_y}^2$.
Since $y$ and $v$ are arbitrary, we have $g\ge Cg_0$ on $B'\cap B^\circ$, and hence $\breve{g}\ge Cg_0$ almost everywhere on $B'$, where $C=1/(C_1C_2)>0$.

Since $\breve{g}$ is locally uniformly positive definite, every form that is $L^2$-integrable with respect to $\breve{g}$ belongs to $L^2_{\mathrm{loc}}$ in the ordinary sense.
Hence the $\dbar$-extension property (see \cite[Lemma 6.8]{Dem82}) is applicable.
As in \cite[Theorem 4.13]{Wat26}, the Nakano semi-positivity of $\breve{g}$ can then be proved by obtaining the optimal $L^2$ estimate on $Y^\circ$, where $\breve{g}|_{Y^\circ}=g$, and extending the resulting solutions across $\Delta$ using the $\dbar$-extension property.
Furthermore, the minimal extension property of $\breve{g}$ follows from Lemma~\ref{Lemma: L2 extension of Hodge sections} by the same argument as in \cite[Proposition 3.6]{Tak24}.
\end{proof}

Griffiths semi-positivity follows from the unique extension of plurisubharmonic functions via the $\limsup$ construction.

\begin{proposition}\label{Proposition: extension of Hodge metric}
The smooth Hodge metric $g$ on $\hdi$ over $Y^\circ$ extends uniquely to a Griffiths semi-positive singular Hermitian metric $\widetilde g$ on $Y^\circ_{q,h}$.
\end{proposition}

\begin{proof}
Fix any $y_0\in Y^\circ_{q,h}$ and take a sufficiently small neighborhood $W\subset Y^\circ_{q,h}$ of $y_0$.
Proposition \ref{Proposition: local uniform positive definiteness} gives a constant $C>0$ such that $\breve g\ge Cg_0$ almost everywhere on $W$, where $g_0$ is any smooth Hermitian metric on $\hdi|_{Y^\circ_{q,h}}$.
Consequently, the dual metrics satisfy $g^*\le C^{-1}g_0^*$ on $W\cap Y^\circ$.
For any $s\in H^0(W,(R^qf_*(K_{X/Y}\otimes\lih))^*)$, put $\varphi_s:=\log|s|_{g^*}^2$. Since $g$ is Nakano semi-positive on $Y^\circ$, the function $\varphi_s$ is plurisubharmonic on $W\cap Y^\circ$.
Furthermore, the preceding estimate gives $\varphi_s\le\log C^{-1}+\log|s|_{g_0^*}^2$, and hence $\varphi_s$ is locally bounded above on $W$.

Thus, the removable singularity theorem for plurisubharmonic functions shows that 
\begin{align*}
  \tl{\varphi}_s(y):=\limsup_{z\in W\cap Y^\circ\!,\, z\to y}\log|s(z)|_{g^*}^2
\end{align*}
is the unique plurisubharmonic extension of $\varphi_s$ to $W$. Define the extended dual metric $\widetilde g^*$ by $|s|_{\widetilde g^*}^2:=e^{\widetilde\varphi_s}$.
Then \cite[Proposition 2.4]{Tak24} implies that $\widetilde g$ is the unique Griffiths semi-positive singular Hermitian metric on $Y^\circ_{q,h}$ extending $g$.

As an alternative proof, the proposition also follows directly from the following argument.
Proposition \ref{Proposition: local uniform positive definiteness} shows that the singular Hermitian extension $\breve g$ of $g$ is Nakano semi-positive and locally uniformly positive definite on $Y^\circ_{q,h}$.
Therefore, the existence of a Griffiths semi-positive extension follows from \cite[Theorem 1.7]{Liu26} and \cite{Mat26}.
\end{proof}

\ref{Main Theorem} follows from Theorem~\ref{Theorem: smooth Nakano semi-positivity} and Propositions~\ref{Proposition: local uniform positive definiteness} and~\ref{Proposition: extension of Hodge metric}.

\begin{remark}
  We answer \cite[Question~1.3-1]{Tak24} as follows: 
  The Griffiths semi-positivity of $\tl{g}$ (or, more precisely, the upper semicontinuity of the $\tl{g}^*\!$-norms of local holomorphic dual sections) implies that $g$ is locally uniformly positive definite up to $\Delta$. 
  We can therefore apply the $\dbar$-extension property \cite[Lemma 6.8]{Dem82} as in the preceding argument and in \cite[Theorem~4.13]{Wat26}. 
  It follows that $\tl{g}$, and more generally any extension $\breve{g}$ of $g$ as a singular Hermitian metric, is Nakano semi-positive on $Y$.
\end{remark}


\medskip
\noindent
{\bf Acknowledgement.} 
The author would like to express his gratitude to Professor Shin-ichi Matsumura for valuable discussions.
The author is supported by Grant-in-Aid for Early-Career Scientists $\sharp$26K16989 from the Japan Society for the Promotion of Science (JSPS).

\medskip
\noindent
{\bf AI disclosure.}
The author used OpenAI's ChatGPT (GPT-6 Sol) as a supplementary aid for
English-language editing, literature searches, and mathematical discussions.
The mathematical ideas, proof strategies, and formulations of the results are the author's own.
The author independently verified all mathematical arguments
and takes full responsibility for the content of this paper.

\begin{thebibliography}{99}
    \bibitem[Ber09]{Ber09}
    B. Berndtsson, \textit{Curvature of vector bundles associated to holomorphic fibrations}. Ann. Math. \textbf{169}, 531-560 (2009).
    \bibitem[BP08]{BP08}
    B. Berndtsson and M. P\u{a}un, \textit{Bergman kernels and the pseudoeffectivity of relative canonical bundles}, Duke Math. J. \textbf{145} (2008), no. 2, 341--378. 
    \bibitem[Blo13]{Blo13}
    Z. B{\l}ocki, \textit{Suita conjecture and the Ohsawa--Takegoshi extension theorem}, Invent. Math. \textbf{193}, 149--158 (2013).
    \bibitem[CGP22]{CGP22}
    J. Cao, H. Guenancia and M. P\u{a}un, \textit{Curvature formula for direct images of twisted twisted relative canonical bundles endowed with a singular metric}, Ann. Fac. Sci. Toulouse Math. \textbf{31} (2022), no 3, 861-905.
    \bibitem[Dem82]{Dem82}
    J.-P. Demailly, \textit{Estimations $L^2$ pour l'op\'{e}rateur $\overline{\partial}$ d'un fibr\'{e} vectoriel holomorphe semi-positif au dessus d'une vari\'{e}t\'{e} K\"{a}hl\'{e}rienne compl\`{e}te}, Ann. Sci. Ec. Norm. Sup. \textbf{15}, (1982), 457-511.
    \bibitem[Dem93]{Dem93}
    J. P. Demailly, \textit{A numerical criterical for very ample line bundles}, J. Differential Geom., \textbf{37} (1993), 323-374.
    \bibitem[Dem-book]{Dem-book}
    J.-P. Demailly, \textit{Complex analytic and differential geometry},  https://www-fourier.ujf-grenoble.fr/demailly/manuscripts/agbook.pdf.
    \bibitem[DPS01]{DPS01}
    J.-P. Demailly, T. Peternell, and M. Schneider, \textit{Pseudo-effective line bundles on compact \kah manifolds}, Internat. J. Math. \textbf{6} (2001), no. 6, 689-741.
    \bibitem[DNWZ23]{DNWZ23}
    F. Deng, J. Ning, Z. Wang and X. Zhou, \textit{Positivity of holomorphic vector bundles in terms of $L^p$-estimates for $\dbar$}, Math. Ann. \textbf{385} (2023), no. 1-2, 575--607. 
    \bibitem[EHN96]{EHN96}
    H. W. Engl, M. Hanke and A. Neubauer, \textit{Regularization of Inverse Problems}, Mathematics and Its Applications, \textbf{375}, Kluwer Academic Publishers, 1996.
    \bibitem[Eno93]{Eno93}
    I. Enoki, \textit{Strong-Lefschetz-type theorem for semi-positive line bundles over compact \kah manifolds}, in Geometry and Global Analysis (Tohoku University, Sendai, 1993), pp. 211-212.
    \bibitem[Fuj04]{Fuj04}
     O. Fujino, \textit{Higher direct images of log canonical divisors}, J. Differential Geom. \textbf{66} (2004), 453--479.
    \bibitem[Fuj78]{Fuj78}
    T. Fujita, \textit{On Kähler fiber spaces over curves}, J. Math. Soc. Japan \textbf{30} (1978), 779--794.
    \bibitem[GR84]{GR84}
    H. Grauert and R. Remmert, \textit{Coherent Analytic Sheaves}, Grund. Math. Wiss., \textbf{265}, Springer, 1984.
    \bibitem[Gri70]{Gri70}
    P. A. Griffiths, \textit{Periods of integrals on algebraic manifolds. III. Some global differential-geometric properties of the period mapping}, Inst. Hautes Études Sci. Publ. Math. \textbf{38} (1970), 125--180.
    \bibitem[GZ15]{GZ15}
    Q. Guan and X. Zhou, \textit{A Proof of Demailly's Strong Openness Conjecture}, Annals of Mathematics, \textbf{182} (2015), pp. 605-616.
    \bibitem[HPS18]{HPS18}
    C. Hacon, M. Popa and C. Schnell, \textit{Algebraic fiber spaces over abelian varieties: around a recent theorem by Cao and P\u{a}un}, in Local and global methods in algebraic geometry, Contemp. Math., \textbf{712}, Amer. Math. Soc., Providence, RI, 2018, pp. 143--195. 
    \bibitem[Ina22]{Ina22}
    T. Inayama, \textit{Nakano positivity of singular Hermitian metrics and vanishing theorems of Demailly-Nadel-Nakano type}, Algebr. Geom. \textbf{9} (1) (2022) 69-92.
    \bibitem[IMW24]{IMW24}
    T. Inayama, S. Matsumura and Y. Watanabe, \textit{Singular Nakano positivity of direct image sheaves of adjoint bundles}, to appear in Anal. PDE; arXiv:2407.11412 (2024). 
    \bibitem[Kat95]{Kat95}
    T. Kato, \textit{Perturbation Theory for Linear Operators}, Classics in Mathematics, Springer, 1995. 
    \bibitem[Kaw81]{Kaw81}
    Y. Kawamata, \textit{Characterization of abelian varieties}, Compositio Math. \textbf{43} (1981), 253--276.
    \bibitem[Kol86a]{Kol86a}
    J. Kollár, \textit{Higher direct images of dualizing sheaves. I}, Ann. of Math. \textbf{123} (1986), 11--42.
    \bibitem[Kol86b]{Kol86b}
    J. Kollár, \textit{Higher direct images of dualizing sheaves. II}, Ann. of Math. \textbf{124} (1986), 171--202.
    \bibitem[Liu26]{Liu26}
    Z. Liu, \textit{Convolution regularization preserves the $L^2$-estimate property for $(1,n)$-forms}, preprint, arXiv:2506.22834v2.
    \bibitem[MY04]{MY04}
    F. Maitani and H. Yamaguchi, \textit{Variation of Bergman metrics on Riemann surfaces}, Math. Ann. \textbf{330} (2004), 477--489.
    \bibitem[Mat22]{Mat22}
    S. Matsumura, \textit{Injectivity theorems with multiplier ideal sheaves for higher direct images under K\"ahler morphisms}, Algebr. Geom. \textbf{9} (2022), no. 2, 122--158.
    \bibitem[Mat26]{Mat26}
    S. Matsumura, \textit{A characterization of Griffiths semipositivity by top-degree optimal $\dbar$-estimates}, preprint, arXiv:2610.25900v1. 
    \bibitem[Mor87]{Mor87}
    A. Moriwaki, \textit{Torsion freeness of higher direct images of canonical bundles}, Math. Ann. \textbf{276} (1987), 385--398.
    \bibitem[MT08]{MT08}
    C. Mourougane and S. Takayama, \textit{Hodge metrics and the curvature of higher direct images}, Ann. Sci. \'{E}cole Norm. Sup. (4) \textbf{41}, 905-924 (2008).
    \bibitem[MT09]{MT09}
    C. Mourougane and S. Takayama, \textit{Extension of twisted Hodge metrics for Kähler morphisms}, J. Differ. Geom. \textbf{83}, 131--161 (2009).
    \bibitem[OT87]{OT87}
    T. Ohsawa and K. Takegoshi, \textit{On the extension of $L^2$ holomorphic functions}, Math. Z. \textbf{195}, 197--204 (1987).
    \bibitem[PT18]{PT18} 
    M. P\u{a}un and S. Takayama, \textit{Positivity of twisted relative pluricanonical bundles and their direct images}, J. Algebraic Geom. \textbf{27} (2018), no. 2, 211--272. 
    \bibitem[PV02]{PV02}
    S. Petermichl and A. Volberg, \textit{Heating of the Ahlfors-Beurling operator: weakly quasiregular maps on the plane are quasiregular}, Duke Math. J. \textbf{112} (2002), no. 2, 281--305.
    \bibitem[QZ26]{QZ26}
    C. Qing and X. Zhou, \textit{Structure theorem, injectivity theorem and vanishing theorem for the cohomology groups of pseudo-effective line bundles}, J. Algebraic Geom. \textbf{35} (2026), no. 2, 269-408.
    \bibitem[SY23]{SY23}
    C. Schnell and R. Yang, \textit{Hodge modules and singular hermitian metrics}, Math. Z. \textbf{303} (2023), no.~2, Paper No.~28, 20 pp.
    \bibitem[Tak24]{Tak24}
    S. Takayama, \textit{Singular Griffiths semi-positivity of higher direct images}, Math. Ann. \textbf{388} (2024) 4343--4353.
    \bibitem[Tak95]{Tak95}
    K. Takegoshi, \textit{Higher direct images of canonical sheaves tensorized with semi-positive vector bundles by proper \kah morphisms}, Math. Ann. \textbf{303} (1995), 389-416.
    \bibitem[Tak97]{Tak97}
    K. Takegoshi, \textit{On cohomology groups of nef line bundles tensorized with multiplier ideal sheaves on compact \kah manifolds}, Osaka J. Math. \textbf{34}(4) (1997) 783-802.
    \bibitem[Vie83]{Vie83}
    E. Viehweg, \textit{Weak positivity and the additivity of the Kodaira dimension for certain fibre spaces}, in Algebraic varieties and analytic varieties (Tokyo, 1981), Adv. Stud. Pure Math. \textbf{1} (1983), 329--353.
    \bibitem[Wat25]{Wat25}
    Y. Watanabe, \textit{Dual Nakano positivity and singular Nakano positivity of direct image sheaves}, Nagoya Math. J. \textbf{257} (2025), 1--24.
    \bibitem[Wat26]{Wat26}
    Y. Watanabe, \textit{Griffiths and Nakano positivity and Nakano-Nadel vanishing theorems for singular Hermitian metrics on complex spaces}, preprint, arXiv:2606.16275v2.
    \bibitem[WZ25]{WZ25}
    Y. Watanabe and Y. Zou, \textit{On the direct image of the adjoint nef and big line bundles}, J. Geom. Anal. \textbf{35} (2025), Paper No. 249, 21pp.
    \bibitem[Wu21]{Wu21}
    X. Wu, \textit{On the hard Lefschetz theorem for pseudoeffective line bundles}, International Journal of Mathematics, \textbf{32}, No. 6 (2021) 2150035, 25pp.
\end{thebibliography}
\end{document}